\documentclass[nohyperref]{article}

\usepackage{microtype}
\usepackage{graphicx}
\usepackage{bm}
\usepackage{algorithm}
\usepackage{algorithmic}
\usepackage{multirow}
\usepackage{amsthm}
\usepackage{amsmath}
\usepackage{subfigure}
\usepackage{amssymb}
\usepackage{booktabs} 

\usepackage{balance}
\usepackage{float}
\usepackage{subeqnarray}
\usepackage{cases}
\usepackage{tabularx}
\usepackage{threeparttable}
\usepackage{makecell}


\usepackage{hyperref}

\usepackage[accepted]{icml2022}

\usepackage[capitalize]{cleveref}

\newtheorem{theorem}{Theorem}[section]

\newtheorem{definition}{Definition}[section]

\newtheorem{lemma}{Lemma}[section]

\newtheorem{remark}{Remark}

\newtheorem{proposition}{Proposition}[section]

\newtheorem{assumption}{Assumption}[section]

\icmltitlerunning{Value Function Based Difference-of-Convex Algorithm for Bilevel Hyperparameter Selection Problems} 

\begin{document}

\twocolumn[
\icmltitle{Value Function Based Difference-of-Convex Algorithm for Bilevel Hyperparameter Selection Problems}




\icmlsetsymbol{equal}{*}

\begin{icmlauthorlist}
	\icmlauthor{Lucy Gao}{waterloo}
	\icmlauthor{Jane J. Ye}{uvic}
	\icmlauthor{Haian Yin}{sustech}
	\icmlauthor{Shangzhi Zeng}{uvic}
	\icmlauthor{Jin Zhang}{sustech,shenzhen}
\end{icmlauthorlist}

\icmlaffiliation{waterloo}{Department of Statistics and Actuarial Science, University of Waterloo, Waterloo, Ontario, Canada}
\icmlaffiliation{uvic}{Department of Mathematics and Statistics, University of Victoria, Victoria, British Columbia, Canada}
\icmlaffiliation{sustech}{Department of Mathematics, SUSTech International Center for Mathematics, Southern University of Science and Technology, Shenzhen, Guangdong, China}
\icmlaffiliation{shenzhen}{National Center for Applied Mathematics Shenzhen, Shenzhen, Guangdong, China}

\icmlcorrespondingauthor{Jin Zhang}{zhangj9@sustech.edu.cn}

\vskip 0.3in
]



\printAffiliationsAndNotice{
	}  

\begin{abstract}
Gradient-based optimization methods for hyperparameter tuning guarantee theoretical convergence to stationary solutions when for fixed upper-level variable values, the lower level of the bilevel program is strongly convex (LLSC) and smooth (LLS). This condition is not satisfied for bilevel programs arising from tuning hyperparameters in many machine learning algorithms. In this work, we develop a sequentially convergent Value Function based Difference-of-Convex Algorithm with inexactness (VF-iDCA). We show that this algorithm achieves stationary solutions without LLSC and LLS assumptions for bilevel programs from a broad class of hyperparameter tuning applications. Our extensive experiments confirm our theoretical findings and show that the proposed VF-iDCA yields superior performance when applied to tune hyperparameters.
\end{abstract}

\section{Introduction}

Virtually all machine learning algorithms require hyperparameter tuning. For example, in regression and classification problems, regularization is often used to induce structured solutions and to control model complexity. The degree of regularization is typically controlled by hyperparameters. Careful tuning of these hyperparameters is critical to the predictive accuracy of the fitted model. 

An effective hyperparameter selection strategy is to choose the hyperparameters that minimize a loss function on a held-out validation set, which amounts to solving a bilevel program (BLP). We consider the following BLP framework:
\begin{equation}\label{original_problem}
	\begin{aligned}
		\min_{ x \in \mathbb{R}^n , \lambda \in \mathbb{R}_+^J} ~~& L(x) \\
		\text{s.t.} ~~& x \in \arg\min_{x' \in  \mathbb{R}^n} \left\{ l(x') + \sum_{i = 1}^{J}\lambda_i P_i(x') \right\},
	\end{aligned}
\end{equation}
where 
$L,l : \mathbb{R}^n \rightarrow \mathbb{R}$ and $P_i : \mathbb{R}^n \rightarrow \mathbb{R}_+,~ i = 1,\ldots,J$ are merely convex (possibly non-smooth) functions.
 In BLP Eq. \eqref{original_problem},
 $\lambda$ is a vector of hyperparameters, the upper-level (UL) problem minimizes the validation error in terms of the hyperparameters, and the lower-level (LL) problem fits a model to training data for a given choice of hyperparameters. Note that in  BLP Eq.  \eqref{original_problem},  we assume that for each $\lambda \in  \mathbb{R}_+^J$,  the LL problem  has solutions.
 
Table \ref{table:model} presents examples of bilevel hyperparameter selection problems of the form \eqref{original_problem} with non-strongly convex and/or non-smooth LL objective functions. In \cref{exten}, we will study a generalization of BLP \cref{original_problem} with LL constraints, which captures the support vector machine classification example in Table \ref{table:model}. 

\begin{table*}[tb]
\caption{Examples of bilevel hyperparameter selection problems of the form \eqref{original_problem} and its generalization. These examples were also explored in \citet{kunapuli2008classification} and \citet{feng2018gradient}. }
\label{table:model}
\vspace{5pt}
\begin{center}
\begin{small}
\resizebox{.99\textwidth}{!}{
\begin{threeparttable}
\begin{tabular}{lccccc}
\toprule
Machine learning algorithm & $x$ & $\lambda$ & $L(x) / l(x)$ & $\sum_{i = 1}^{J}\lambda_i P_i(x)$ & LL Constraints \\
\midrule
elastic net
& $\bm{\beta}$ & $\lambda_1$, $\lambda_2$
& $\frac12 \sum_{i \in I_{\text{val}}/i\in I_{\text{tr}}} | b_i - \bm{\beta}^\top \mathbf{a}_i |^2$ 
& $\lambda_1\|\bm{\beta}\|_1 + \frac{\lambda_2}2 \|\bm{\beta}\|_2^2$
& - \\ \specialrule{0em}{3pt}{3pt}
sparse group lasso
& $\bm{\beta}$ & $\lambda\in \mathbb{R}^{M+1}_+$
& $\frac12 \sum_{i \in I_{\text{val}}/i\in I_{\text{tr}}} | b_i - \bm{\beta}^\top \mathbf{a}_i |^2$ 
& $\sum_{m=1}^M \lambda_m \|\bm{\beta}^{(m)}\|_2 + \lambda_{M+1}\|\bm{\beta}\|_1$
& - \\ \specialrule{0em}{3pt}{3pt}
low-rank matrix completion
& $\bm{\theta},\ \bm{\beta},\ \Gamma$ & $\lambda\in \mathbb{R}^{2G+1}_+$
& $\sum_{(i, j) \in \Omega_{\text{val}} / (i, j) \in \Omega_{\text{tr}}} \frac12|M_{ij} - \mathbf{x}_i \bm{\theta} - \mathbf{z}_j \bm{\beta} - \Gamma_{ij} |^2$
& $\lambda_0 \|\Gamma\|_* + \sum_{g=1}^G \lambda_g \|\bm{\theta}^{(g)}\|_2 + \sum_{g=1}^G \lambda_{g+G} \|\bm{\beta}^{(g)}\|_2$
& - \\ \specialrule{0em}{3pt}{3pt}
support vector machine
& $\mathbf{w},\ c$ & $\lambda$, $\bar{\mathbf{w}}$
& $\sum_{j\in I_{\text{val}} / j \in I_{\text{tr}} }\max( 1 - b_j(\mathbf{w}^\top \mathbf{a}_j - c),0)$ 
& $\frac{\lambda}{2}\|\mathbf{w}\|^2 $
& $-\bar{\mathbf{w}} \le \mathbf{w} \le \bar{\mathbf{w}}$\\
\bottomrule
\end{tabular}
\begin{tablenotes}
	\item[] The detailed descriptions of these problems are presented in the Supplemental Material.
\end{tablenotes}
\end{threeparttable}}
\end{small}
\end{center}
\vskip -0.2in
\end{table*}

\subsection{Related Work}

The simplest hyperparameter selection approach is brute-force grid search, which is seldom applied when there are more than two hyperparameters due to its unfavourable computational complexity. The current gold standard in the presence of more than two hyperparameters is Bayesian optimization, which models the conditional probability of the performance on some metric given hyperparameters and a dataset. However, these gradient-free methods often scale poorly beyond 10 to 20 hyperparameters. 

Gradient-based optimization methods (GM) can handle up to a few hundred hyperparameters. Existing gradient-based methods can be roughly divided into two categories based on the strategy used for calculating the hypergradient (the gradient
of the UL objective). Implicit gradient methods (IGM), also known as implicit differentiation~\cite{pedregosa2016hyperparameter,rajeswaran2019meta,lorraine2020optimizing},  replace the LL sub-problem with an implicit equation. These IGM apply the implicit function theorem (IFT) to the optimality conditions of the LL
problem, and hence derive their hypergradients by solving a linear system.
As it involves computing a Hessian matrix and its inverse, in practice, the conjugate gradient method (CG) or Neumann method are used for fast inverse computation.
Explicit gradient methods (EGM) replace the LL problem with a gradient-descent-type dynamics iteration~\cite{franceschi2017forward,franceschi2018bilevel,finn2017model,nichol2018first,lorraine2018stochastic}. 
Specifically, \citet{franceschi2017forward} and \citet{franceschi2018bilevel} first calculate gradient representations of the LL objective and then perform either reverse or forward gradient computations (called Reverse Hyper-Gradient (RHG) and Forward Hyper-Gradient (FHG)) for the UL objective. 
In order to reduce the amount of computation,~\citet{shaban2019truncated} proposes truncated reverse hypergradient to truncate the gradient trajectory, and \citet{liu2018darts} uses the difference of vectors to approximate the gradient.
Recently, \citet{ji2020convergence} and \citet{ji2021bilevel} analyze the convergence rate and complexity of IGM and EGM, respectively. \citet{ji2022will} analyze the effect of inner loop step for IGM and EGM.
Besides, some work \cite{liu2020generic, liu2021value, sow2022constrained, liu2022general} study the case where LL problem is convex but not strongly convex; \citet{Liu2021TowardsGB} proposed an initialization auxiliary method for the setting where the LL problem is generally nonconvex. 

\subsection{Our Motivations and Contributions}

Some theoretical progress has been made in both gradient-free and gradient-based methods. However, existing GM require restrictive assumptions to ensure the convergence towards stationarity: both IGM and EGM require the LL strong convexity (LLSC) and smoothness (LLS), which is often violated in applications. For IGM, the subsequential convergence of stationary solutions can be found in~\citet{pedregosa2016hyperparameter}.
Without the LLS assumption, some IGM variants study the computation of hypergradients, where the convergence analysis of stationarity of such methods is still lacking; see Implicit Gradient-based Joint Optimization (IGJO) in \citet{feng2018gradient} and Implicit Forward Differentiation Method (IFDM) in \citet{bertrand2020implicit} and   \citet{bertrand2021implicit}.
For EGM,
when the LL objective is strongly convex uniformly for $\lambda$ and its gradient w.r.t. $x$ is uniformly Lipschitz continuous, the LL gradient iterates converge linearly in a uniform rate. This further leads to a uniformly linear convergence of the hypergradient of approximated UL objective to the true hypergradient.
The subsequential convergence  towards stationary solutions for EGM thus follows easily; see ~\cite{grazzi2020iteration}.
Focusing on stationarity convergence, under LLSC and LLS, ~\citet{Ji2020} and \citet{Ji2021} give a comprehensive characterization on the nonasymptotic convergence rate and complexity bounds for both IGM and EGM.

Many applications involve BLP where \emph{LLSC and LLS are violated}; see the hyperparameter selection applications in BLP \cref{original_problem} and Table \ref{table:model}. This motivates us to investigate algorithms with provable convergence towards stationary solutions \emph{without} assuming LLSC and LLS. 
\vskip -0.1 in

\begin{table}[th]\vspace{-5pt}
\caption{Comparing theoretical results among BLP algorithms. }
\label{table:theoretical_comparison}
\vspace{-5pt}
\begin{center}
\begin{small}
\resizebox{.49\textwidth}{!}{
\begin{tabular}{lc cc cc}
\toprule
Category & Methods & LLSC & LLS & Conv. Type & Sol. Quality \\
\midrule
EGM & FHG/RHG & w/ & w/ & subsequential & stationary \\
\specialrule{0pt}{3pt}{3pt}
\multirow{2}{*}{IGM} & CG/Neumann & w/ & w/ & subsequential & stationary \\
\specialrule{0pt}{1pt}{1pt}
& IGJO/IFDM & w/ & w/o & - & - \\
\specialrule{.5pt}{3pt}{3pt}
DCA & VF-iDCA & w/o & w/o & sequential & stationary \\
\bottomrule
\end{tabular}
}
\end{small}
\end{center}
\vskip -0.1in
\end{table}

A striking feature of our study is the
decoupling of  the hyperparameters from the regularization term, which results in an equivalent BLP reformulation with fully convex LL problem. 
Critically, since the reformulated LL problem is fully convex, the value function (VF) of the LL problem is convex. Subsequently, we can further reformulate BLP \cref{original_problem} as a DC program. By using the main idea of DC algorithm (DCA) which linearizes the concave part of the DC structure, we propose a sequential convex programming scheme with subproblem inexactness, named VF-iDCA. For a wide range of applications, we justify the sequential convergence towards stationary solutions, which is an unusual finding in BLP optimization.  
Recently, a DC-type algorithm, namely, iP-DCA for solving bilevel support vector classification problem was proposed in~\citet{DCbilevel}. In comparison with iP-DCA, the new hyperparameter decoupling mechanism of VF-iDCA allows us to efficiently handle more hyperparameter selection BLPs where
the LL  problems are equipped with complex regularization terms.
Our main contributions are as follows: 
\vspace{-5pt}
	\begin{itemize}
         \item By decoupling hyperparameters from the regularization, based on the value function approach, VF-iDCA constructs a series of tractable convex subproblems.

		\item  To our best knowledge, we establish the first strict sequential convergence towards stationary solutions for an important class of hyperparameter selection applications without assuming LLSC and LLS for LL tasks. 
		\item  We conduct experiments to verify our theoretical findings and evaluate VF-iDCA
on various hyperparameter selection BLPs. Our results compare favourably to  existing gradient-free and gradient-based methods. 
	\end{itemize}
\vspace{-5pt}
In Table \ref{table:theoretical_comparison}, we summarize our contributions to the theory of convergence of BLP algorithms by comparing our results to some existing results in the literature.

\section{Value Function Based Difference-of-Convex Algorithm with Inexactness}

\subsection{Hyperparameter Decoupling and Fully Convex LL}
The LL problem in \cref{original_problem} is a penalized problem with regularization term involving coupled UL and LL variables:
\begin{equation}
	~~~~~~\min_{x'}\  l(x') + \sum_{i = 1}^{J}\lambda_i P_i(x'). \label{convexop1}
\end{equation}
We decouple the hyperparameter variables $\lambda$ from the regularization term by introducing a new variable $r$, which yields the following constrained optimization problem:
\begin{equation}
	~~~~~~\min_{x' }\ l(x')~~~\text{s.t.}~ P_i(x') \le r_i,~ i=1,\ldots,J. \label{convexop2}
\end{equation}
Denote their solution sets by $ S_p(\lambda)$ and $S_c(r)$, respectively.
Thanks to the LL convexity w.r.t. $x$, there is a one-to-one relationship between regularized problem \cref{convexop1} and constrained problem \cref{convexop2}. That is, for any $\lambda \ge 0$, there is $r \ge 0$ such that $S_p(\lambda) = S_c(r)$ and vice versa. This suggests working with the following BLP:
\begin{equation}\label{reform_problem}
	\begin{aligned}
		\min_{x, r \in \mathbb{R}_+^J} ~~& L(x) \\
		\text{s.t.} ~~& x \in \arg\min_{x'} \{l(x')~\text{s.t.}~ P_i(x') \le r_i, i=1,\ldots,J \}.
	\end{aligned}
\end{equation}
In the section \ref{sec:alg}, we will use the fact that the LL of \cref{reform_problem} is \emph{fully convex w.r.t the joint variable $(x, r)$} to rewrite \cref{reform_problem} as a single-level problem and design an algorithm to solve it. 

However, we must ask: is \cref{reform_problem} equivalent to \cref{original_problem}? It is relatively straightforward to show that their global solutions coincide, using the equivalence between their LL problems. However, just because their global solutions coincide, does not necessarily mean their local solutions coincide. If we are to work with \cref{reform_problem} in place of \cref{original_problem}, then it is critical that the local solutions to \cref{reform_problem} are local solutions to \cref{original_problem}, because \cref{reform_problem} is a non-convex optimization problem. Theoretical analysis of the relationship between the local solutions is substantially more challenging than theoretical analysis of the relationship between the global solutions. Fortunately, we are able to overcome this challenge in Proposition \ref{local_min1} of Section \ref{sec:solution}, where we show that local minima of \cref{reform_problem} are also local minima of \cref{original_problem}.

\subsection{Single-level DC Reformulation and Algorithm}
\label{sec:alg}

We now consider solving \eqref{reform_problem}.
The value function of the LL problem governed by $r$ is denoted by $$v(r) := \min \left\{l(x)~\text{s.t.}~ P_i(x) \le r_i,~ i=1,\ldots,J \right\}.$$ 
Throughout the paper, we assume that the following condition holds.
\begin{assumption}\label{assum1} For each $r$ belonging to any open subset of $ \mathbb{R}_+^J$, 
 $\mathcal{F}(r):=\{x ~\text{s.t.}~ P_i(x)\leq r_i, i=1,\dots, J\}\not =\emptyset$ and $l(x)$ is bounded below. 
\end{assumption}

Thanks to full convexity, $v(r)$ induced by a partial minimization is convex and locally Lipschitz continuous around every point in $\mathbb{R}_{>0}^J:= \{r \in \mathbb{R}^J | r_j > 0\}$ (see e.g. Lemma 3 in \citealt{DCbilevel}). 
Using the value function, we now reformulate BLP  \eqref{reform_problem} as the following  DC program:
\begin{equation}\label{p-singlelevel_problem}
	\begin{aligned}
		\min_{x, r \in \mathbb{R}_+^J} ~~& L(x) \\
		\text{s.t.} ~~~& l(x) - v(r) \le 0, 
		 P_i(x) \le r_i,~ i=1,\ldots,J.
	\end{aligned}
\end{equation}

We next propose VF-iDCA to 
solve BLP \eqref{reform_problem}. By using the main idea of DCA which linearizes the concave part of the
DC structure, we propose a sequential convex programming scheme as follows.
Given a current iteration $({x}^k,r^k)$ for each $k$, solving the LL problem parameterized by $r^k$
\begin{equation}\label{p-LL-iter}
	\min_{x} ~~ l(x)~\text{s.t.}~ P_i(x) \le r^k_i,~ i=1,\ldots,J,
\end{equation}
leads to a solution $\tilde{x}^k \in S_c(r^k)$ and a corresponding Karush-Kuhn-Tucker (KKT) multiplier $\gamma^k \in \mathcal{M}(\tilde{x}^k,r^k)$, where $\mathcal{M}(x,r)$ denotes the set of KKT multipliers of LL problem
\[
\begin{aligned}
\mathcal{M}(x,r):= & \Big\{ \gamma \in  \mathbb{R}^J_+\mid 0\in \partial l(x) +\sum_{i=1}^J \gamma_i \partial P_i(x),  \\
&~ \gamma_i (P_i(x)-r_i)=0, i=1,\dots, J \Big\}.
\end{aligned}
\]

Then by sensitivity analysis, we have  $- \gamma^k \in \partial v(r^k)$; see Proposition \ref{Prop3.3}.
Compute $z^{k+1} := (x^{k+1},r^{k+1})$ as an approximate minimizer of the strongly convex subproblem			
\begin{equation} \label{p-DCA2_subproblem}
	\begin{aligned}
		\min_{x, r \in \mathbb{R}_+^J} ~~ \phi_k(x,r) := &L(x) +  \frac{\rho}{2} \| z - z^k \|^2 \\
		 & +\alpha_k \max_{ i=1,\ldots,J}\{0, V_k(x,r), P_i(x) - r_i\},
	\end{aligned}
\end{equation}
where $\rho>0$, and $\alpha_k$ represents the adaptive penalty parameter,   $z:= (x,r)$, $z^k := (x^k,r^k)$ and $$V_k(x,r):= l(x) - l(\tilde{x}^k) +\langle \gamma^k, r-r^k\rangle.$$ Denoting $\Sigma:= \mathbb{R}^n \times \mathbb{R}_+^J $, 
we introduce an inexact condition for choosing $z^{k+1}$:
\begin{equation}\label{p-inexact2}
	\mathrm{dist}(0, \partial \phi_k(z^{k+1}) + \mathcal{N}_\Sigma(z^{k+1})) \le \frac{\sqrt{2}}{2} \rho\|z^k - z^{k-1}\|,
\end{equation}
where $\mathrm{dist}(x, \Omega) $ is the distance from $x$ to $\Omega$ and  $ \mathcal{N}_\Sigma$ denotes the normal cone to $\Sigma$.
Using above constructions and
letting \begin{equation} 
t^{k+1} = \max_{ i=1,\ldots,J}\{0, V_k(x^{k+1},r^{k+1}), P_i(x^{k+1}) - r^{k+1}_i\},\label{deft}
\end{equation} 
we are ready to present  VF-iDCA in \cref{p-ipDCA}.
\begin{algorithm}[h]
	\caption{VF-iDCA}\label{p-ipDCA}
	\begin{algorithmic}[1]
		\STATE Take an initial point $(x^0,r^0)\in \mathbb{R}^n\times  \mathbb{R}_+^J$; $c_\alpha, \delta_\alpha > 0$; an initial penalty parameter $\alpha_0 > 0$; tolerance $tol > 0$.
		\FOR{$k=0,1,\ldots$}
				\STATE Solve LL problem \cref{p-LL-iter}. Find
			$\tilde{x}^k \in S_c(r^k)$ and a KKT multiplier $\gamma^k$.		
				\STATE Solve problem \cref{p-DCA2_subproblem} up to tolerance in \eqref{p-inexact2}. Find an approximate solution $z^{k+1}=(x^{k+1},r^{k+1})$.
				\STATE Stopping test: Stop if $\max\{ \|z^{k+1} - z^k\|,t^{k+1} \} < tol$.
				\STATE Adaptive penalty parameter update:
				Set
				\begin{equation} \label{adaptpen}
					\alpha_{k+1} = \left\{
					\begin{aligned}
						&\alpha_k + \delta_\alpha, &&\text{if}~\max \Big\{\alpha_k, \frac{1}{t^{k+1}} \Big\} < \frac{c_\alpha}{\Delta^{k+1}}, \\
						&\alpha_k, &&\text{otherwise},
					\end{aligned}\right.
				\end{equation}
			 for $t^{k+1}$ defined in \eqref{deft} and $\Delta^{k+1}:=\|z^{k+1} - z^k\|$. 
		\ENDFOR
	\end{algorithmic}
\end{algorithm}

\begin{remark}
When $l(x)$ is differentiable with a Lipschitz continuous gradient, we can derive a linearized version of VF-iDCA  where we not only
linearize the concave part but also the convex smooth part of the DC structure. 
\end{remark}

\section{Theoretical Investigations}\label{sec:theory}

With the purpose of deriving a provably convergent algorithm for BLP \cref{original_problem} without LLSC and LLS, we involve two signature features in our algorithmic design, i.e., hyperparameter decoupling and
single-level DC reformulation. 
In this section, we provide convergence analysis of the proposed VF-iDCA towards ``good" quality solutions, i.e. towards KKT stationarity solutions, defined as follows.
\begin{definition} A feasible point $(\bar x, \bar r)$ of problem \cref{p-singlelevel_problem} is a KKT stationary solution  if there exist multipliers $\eta\geq 0, \lambda \in \mathbb{R}_+^J$ such that
		\begin{equation*}
		\begin{aligned}
				& 0 \in \partial L(\bar x)  + \eta \partial l(\bar{x}) + \sum_{i=1}^J \lambda_i \partial P_i(\bar x), \\
				&  \lambda\in  -\eta \partial v(\bar r) +{\cal N}_{\mathbb{R}_+^J}(\bar r),  \lambda_i \left(  P_i(\bar x) - \bar{r}_i \right) = 0, i=1,\dots, J.
		\end{aligned}
	\end{equation*}

\end{definition}

\subsection{Solutions Recovery for Hyperparameter  Decoupling}
\label{sec:solution}
Before analyzing the convergence properties of Algorithm \ref{p-ipDCA}, we first need to justify the validity of the hyperparameter decoupling, i.e. we need to show that BLP \cref{original_problem} is equivalent to BLP \cref{reform_problem}.

The following lemmas show a one-to-one relationship between the LL problems \cref{convexop1} and \cref{convexop2}.
\begin{lemma}\label{p-inclusion1}
	Let $x \in S_p(\lambda)$ with  $\lambda \in \mathbb{R}_+^J$.   Then $x \in S_c(r)$ where  $r = P(x) \in \mathbb{R}_+^J$. 
\end{lemma}

\begin{lemma}\label{p-inclusion2}
	For any  $x \in S_c(r)$ with $r \in \mathbb{R}_+^J$, suppose 
	$\lambda\in \mathcal{M}(x,r)$. Then $ x \in S_{p}(\lambda)$. 
\end{lemma}
Note that the nonemptyness of the multiplier set $\mathcal{M}(x,r)$ can be acheived by imposing certain constraint qualifications such as the Slater condition.

The next theorem shows that the $x$ components of globally optimal solutions for BLPs \cref{original_problem} and \cref{reform_problem} coincide.
We denote by  $P(x):=(P_1(x),\dots, P_J(x))$. 
\begin{theorem}[Global solutions recovery]\label{p-globalequi}
If $(\bar{x}, \bar{r}) \in \mathbb{R}^n \times \mathbb{R}_+^J$  is a global optimal solution of problem \cref{reform_problem} and $\bar \lambda \in \mathcal{M}(\bar x,\bar r)$,
then $(\bar{x}, \bar \lambda)$ is a global optimal solution of BLP \cref{original_problem}. 
Conversely, suppose there exists a dense subset $D$ of $\mathbb{R}^n \times \mathbb{R}_+^J$ such that $\mathcal{M}(x,r)\not =\emptyset$ for all $(x,r) \in D$. Let $(\bar{x}, \bar \lambda)$ be a global optimal solution of problem  \cref{original_problem}. Then 
$(\bar{x}, \bar r) $ with $\bar r:=P(\bar{x})$ is a global  optimal solution of BLP \cref{reform_problem}.
\end{theorem}

Theorem \ref{p-globalequi} says that the global solutions of the fully convex BLP \cref{reform_problem} are the same as the global solutions of the original BLP \cref{original_problem}.

We now ask an arguably more important question: are the \emph{local} solutions of BLP \cref{reform_problem} the same as the local solutions of BLP \cref{original_problem}? First, we show in the following proposition that any local optimal solution of BLP \cref{reform_problem} must be locally optimal for BLP \cref{original_problem}. 
\begin{proposition}[Local solution recovery-I]\label{local_min1}
	If $(\bar{x},\bar{r})$ is a local solution of BLP \cref{reform_problem} with $\bar{r} = P(\bar{x})$ and $\bar \lambda \in \mathcal{M}(\bar x,\bar r)$,
	then $(\bar{x},\bar{\lambda})$ is a local solution of BLP \cref{original_problem}.
\end{proposition}
This is the key result that justifies our choice in this paper to solve the tractable fully convex BLP \cref{reform_problem} in place of the substantially more difficult BLP \cref{original_problem}. 

The reverse direction is less practically relevant, because we are unlikely to attempt to solve BLP \cref{reform_problem} by solving the more difficult BLP \cref{original_problem}. Nevertheless, the following result provides conditions under which local solutions to BLP \cref{original_problem} are local solutions to BLP \cref{reform_problem}. 
\begin{proposition}[Local solution recovery-II]\label{local_min2}
	Suppose that there exists a dense subset $D$ of $\mathbb{R}^n \times \mathbb{R}_+^J$ such that $\mathcal{M}(x,r)\not =\emptyset$ for all $(x,r) \in D$. Let $(\bar{x}, \bar \lambda)$ be a local optimal solution of BLP \cref{original_problem},  locally w.r.t. $x$ and globally w.r.t. $\lambda$. Then 
	$(\bar{x}, \bar r) $ with $\bar r:=P(\bar{x})$ is an optimal solution of BLP \cref{reform_problem} locally w.r.t. $x$ and globally w.r.t. $r$.
	
	Furthermore, let $(\bar{x},\bar \lambda)$ be a locally optimal solution of  BLP \cref{original_problem} and $\bar r:=P(\bar{x})$. If there exist a bounded set $\Lambda$ and $\epsilon>0$ such that $\Lambda \cap \mathcal{M}(x,r) \neq \varnothing$ for any $(x,r)\in D \cap \mathbb{B}_{\epsilon}(\bar x,\bar{r})$,  $\mathcal{M}(\bar x,\bar r)$ is a singleton, then $(\bar{x}, \bar r)$ must be locally optimal for BLP \cref{reform_problem}.
\end{proposition}
\begin{remark} If both $l, P$ are differentiable, then the existence of dense set and the uniqueness of the multipliers can be simply implied by the  linear independence of the gradient vectors
	$\{\nabla  P_i(\bar{x})\}_{i=1}^J.$
\end{remark}

\subsection{Sequential Convergence Towards Stationary Solutions}
\begin{assumption}\label{p-assum} Throughout this section, we assume the following conditions.
	\begin{itemize}
		\item[(a)] For all the $r^k$  generated by the algorithm, the lower level solution set $S_c(r^k)$ is nonempty.
		\item[(b)] For all the $r^k$ and  $\tilde{x}^k \in S_c(r^k)$ generated by the algorithm, $\mathcal{M}(\tilde{x}^k, r^k)\not = \emptyset$. 
	\end{itemize}
\end{assumption}

Under  \cref{assum1},  by Lemma 3 in \citet{DCbilevel}, the value function is locally Lipschitz continuous in $\mathbb{R}_{>0}^J$.
By Theorem 3 in \citet{DCbilevel}, 
we have a characterization for the subdifferential of $v(r)$ in the sense of convex analysis~\cite{rockafellar2015convex}.
\begin{proposition}\label{Prop3.3}
	For any $r \in \mathbb{R}_+^J$, $x\in S_c(r)$, we have that\vspace{-2pt}
	\begin{equation*}
			\partial v(r) \supseteq  -\mathcal{M}(x, r). \vspace{-2pt}
	\end{equation*} 
	In addition, if  the Slater condition holds for the LL problem, then the equality holds.
\end{proposition}

The subsequential convergence of VF-iDCA towards stationary solutions follows from Theorem 1 in \citet{DCbilevel}.

\begin{theorem}\label{p-thm_limit} Suppose that  $L(x)$ is bounded below and the sequences $\{z^k \}$ and $\{\alpha_k\}$  generated by the VF-iDCA are bounded.  Then any accumulation point 
 $(\bar{x},\bar{r})$  of
 $\{z^k\}$ is a KKT stationary solution of problem \cref{p-singlelevel_problem} provided that   $\bar{r} \in \mathbb{R}_{>0}^J$.
\end{theorem}

Next, as the main contribution of this part, we show that the subsequential convergence  can be further enhanced to sequential convergence under the Kurdyka-\L{}ojasiewicz property \citep{attouch2009convergence,attouch2010proximal,attouch2013convergence,bolte2014proximal}. Let $\eta \in [0, +\infty]$ and $\Phi_\eta$ denote the class of all concave and continuous functions $\varphi: [0,\eta) \rightarrow [0, +\infty)$ satisfying the conditions: (a) $\varphi(0) = 0$, (b) $\varphi$ is $C^1$ on $(0,\eta)$ and continuous at $0$, (c) $\varphi'(s) > 0$ for all $s \in (0,\eta)$.
\begin{definition}[Kurdyka-\L{}ojasiewicz property]
	Let $\sigma : \mathbb{R}^d \rightarrow (-\infty, +\infty]$	be proper and lower semicontinuous.
	The function $\sigma$ is said to have the Kurdyka-\L{}ojasiewicz (KL) property
	at $\bar{y} \in \mathrm{dom}\, \partial \sigma := \{ y \in \mathbb{R}^d\ ~|~  \partial \sigma(y) \neq \varnothing \}	$,	
	if there exist $\eta \in (0,+\infty]$, a neighborhood $\mathcal{Y}$ of $\bar{y}$ and a function $\varphi \in \Phi_\eta$, such that for all
	$y \in \mathcal{Y} \cap \{ y \in \mathbb{R}^d\ ~|~  \sigma(\bar{y}) < \sigma(y) <\sigma(\bar{y}) +\eta \}$,
	the following inequality holds\vspace{-2pt}
	\[
	\varphi'(\sigma(y) - \sigma(\bar{y})) \mathrm{dist}(0, \partial \sigma(y)) \ge 1.\vspace{-2pt}
	\]
	If $\sigma$ satisfy the KL property at each point of $\mathrm{dom}\, \partial \sigma$ then $\sigma$ is called a KL function.
\end{definition}

In addition, given the KL property around any points in a compact set, the uniformized KL property holds; see Lemma 6 in \citet{bolte2014proximal}.
\begin{lemma}[Uniformized KL property] \label{p-uniformKL}
	Given a compact set $C$ and a proper and lower semicontinuous function $\sigma : \mathbb{R}^d \rightarrow (-\infty, +\infty]$, suppose
	that $\sigma$ is constant on $C$ and satisfies the KL property at each point of $C$. Then,
	there exist $\epsilon, \eta$ and $\varphi \in \Phi_\eta$ such that for all $\bar{y} \in C$ and	$y \in\{ y \in \mathbb{R}^d\ ~|~  \mathrm{dist}(y, C) < \epsilon, ~~ \sigma(\bar{y}) < \sigma(y) <\sigma(\bar{y}) +\eta \}$, it holds \vspace{-2pt}
	\[
	\varphi'(\sigma(y) - \sigma(\bar{y})) \mathrm{dist}(0, \partial \sigma(y)) \ge 1.
	\]
\end{lemma}\vspace{-2pt}

Inspired by \citet{liu2019refined}, we define the following merit function for convergence analysis:\vspace{-2pt}
\[
\begin{aligned}
	E_{\alpha}(z,z_0, \gamma):=& L(x) + \frac{\rho}{4} \|z-z_0\|^2 + \delta_{\Sigma}(z)\\
	& \hspace{-55pt} +\alpha \max_{ i=1,\ldots,J}\{0, l(x) + \langle \gamma, r \rangle + v^*(-\gamma),P_i(x) - {r}_i\},
\end{aligned}
\] 
where $z:= (x,r)$, $z_0 := (x_0,r_0)$, $ \delta_{\Sigma}(z)$ is 
the indicator function of $\Sigma$ and $v^*$ denotes the conjugate function of $v$, that is, $v^*(\xi) := \sup_r \{\langle \xi,r\rangle - v(r) \}$.

Let $\{(x^k, r^k)\}$ be iterates generated by VF-iDCA. Before we can analyze the sequential convergence of $\{(x^k, r^k)\}$, we need the sufficient decrease property and the relative error condition of the merit function $E_{\alpha}$, which are summarized in the following lemma.
\begin{lemma}\label{p-suff_decreasenew}
	Let $\{(x^k, r^k)\}$ be iterates generated by VF-iDCA, then
	$z^k := (x^k, r^k)$ satisfies 
	\begin{equation} 
		\begin{aligned}
				&E_{\alpha_k}(z^k,z^{k-1},\gamma^{k-1}) \\ \geq\, &E_{\alpha_k}(z^{k+1},z^{k},\gamma^{k}) + \frac{\rho}{4} \|z^{k+1} - z^k\|^2, \label{p-suff_decrease_eq2} \quad \text{and}
		\end{aligned}
		\end{equation}
\begin{equation}\label{p-re_err_eq}
		\begin{aligned} 
			&\mathrm{dist} \left( 0, \partial E_{\alpha_k}(z^{k+1},z^{k},\gamma^{k}) \right) \\ \le \, & \frac{\sqrt{2}}{2} \rho\|z^k - z^{k-1}\| +  (\alpha_k+\rho) \|z^{k+1} - z^k\|.
		\end{aligned}
	\end{equation}
\end{lemma}

\begin{theorem}\label{p-con_KL_alg1}
	Suppose that  $L(x)$ is bounded below and the sequences $\{z^k \}$ and $\{\alpha_k\}$ generated by VF-iDCA are bounded,  the merit function $E_\alpha$ is a KL function and there exists $\delta > 0$ such that $r_i^k \ge \delta$ for all $k$ and $i=1,\dots, J$.
	Then the sequence $\{z^k\}$ converges to a KKT stationary solution of problem \cref{p-singlelevel_problem}.
\end{theorem}

\subsection{Further Justification of Validity in Applications}

To further justify the validity of our sequential convergence theory in real-world applications, we next discuss the conditions, especially the KL property imposed on the merit function $E_\alpha$ in \cref{p-con_KL_alg1}.
First, part (b) of \cref{p-assum} holds automatically, for example if the Slater condition holds, i.e., there exists $x_0$ such that $P_i(x_0)<r^k_i$ for all $i=1,\dots, J$.
Thus, a large class of functions automatically satisfies the KL property, see  e.g.  \citet{bolte2010characterizations}, \citet{attouch2010proximal}, \citet{attouch2013convergence}, and \citet{bolte2014proximal}. As shown in \citet{bolte2007lojasiewicz} and \citet{bolte2007clarke}, any semi-algebraic function meets  the KL property.
\begin{definition}[Semi-algebraic sets and functions]
	A subset $S \subseteq \mathbb{R}^d$ is called a real semi-algebraic set if there exists a finite number
	of real polynomial functions $g_{ij}, h_{ij}: \mathbb{R}^d \rightarrow \mathbb{R}$ such that
	\[
	S = \bigcup_{j = 1}^p\bigcap_{i = 1}^q\{y \in \mathbb{R}^d ~|~ g_{ij}(y) = 0, ~ h_{ij}(y) < 0 \}.
	\]
	A function $h : \mathbb{R}^d \rightarrow (-\infty, +\infty]$ is called semi-algebraic if its graph
	$\{(y,s) \in \mathbb{R}^{d+1} ~|~ h(y) = s\}$
	is a semi-algebraic subset of $\mathbb{R}^{d+1}$.
\end{definition}

\begin{lemma}
	Let $\sigma : \mathbb{R}^d \rightarrow (-\infty, +\infty]$	be proper and lower semicontinuous. If $\sigma$ is semi-algebraic then it satisfies the KL property at any point of its domain.
\end{lemma}

The class of semi-algebraic functions is closed under various operations, see, e.g., \cite{attouch2010proximal,attouch2013convergence,bolte2014proximal}. In particular, the indicator functions of semi-algebraic sets, finite sum and product of semi-algebraic functions, composition of semi-algebraic functions and partial minimization of semi-algebraic function over semi-algebraic set
are all semi-algebraic functions. Applying this to BLP \cref{original_problem}, when the LL functions $l(x)$ and $P(x)$ are both semi-algebraic, $v(r)$ is a semi-algebraic function. Moreover, if $L(x)$ is a semi-algebraic function, the merit function $E_\alpha$ is also semi-algebraic and hence a KL function. 

\begin{theorem}
	Assume that $L(x)$, $l(x)$ and $P(x)$ are semi-algebraic functions. Suppose that $\{z^k := (x^k, r^k)\}$ and $\{\alpha_k\}$ generated by VF-iDCA are bounded,  $L(x)$ is bounded below and there exists $\delta > 0$ such that $r^k_i \ge \delta$ for all $k$ and $i=1,\dots, J$. Then $\{z^k\}$ converges to a KKT point of problem \cref{p-singlelevel_problem}.
\end{theorem}

A wide range of functions appearing in applications are semi-algebraic, see, e.g., \citet{attouch2010proximal}, \citet{attouch2013convergence}, and \citet{bolte2014proximal}. In particular, all the functions involved in our applications in the experiments section are semi-algebraic. 

\section{Extension to General Setting}\label{exten}
Our valued function based DC algorithm as well as its convergence analysis can be straightforwardly extended to BLP in a more general setting with LL constraints.
\begin{equation}\label{p-general_problem}
	\begin{aligned}
		\min_{x \in X, u \in U, \lambda \in \mathbb{R}_+^J} ~~& L(x,u) \\
		\text{s.t.} ~~ x \in \underset{x' \in X}{\mathrm{argmin}} &\bigg\{ l(x',u) + \sum_{i = 1}^{J}\lambda_i P_i(x',u), \\[-10pt]
		   & \hspace{80pt} 
		\text{s.t.} ~~ g(x', u) \le 0 \bigg\},
	\end{aligned}
\end{equation}
where $X \subseteq \mathbb{R}^n$, $U \subseteq \mathbb{R}^d$ are closed convex sets,  $g_i(i=1,\dots,m),L,l , P_i:\mathbb{R}^n \times \mathbb{R}^d \rightarrow \mathbb{R}_+(i=1,\dots,J )$ 
are  convex functions defining on an open convex set containing  $X\times U$ and $g(x,u) := (g_1(x,u),\ldots,g_J(x,u))$.
We provide a detailed description of VF-iDCA for BLP \cref{p-general_problem} in the Supplemental Material.

\section{Experiments}

In this section we will compare the performance of our VF-iDCA to widely-used competitors in hyperparameter optimization for a number of machine learning algorithms. Detailed descriptions of these problems are presented in the Supplemental Material. All algorithms were implemented in Python and the software package used for reproduce our experiments is available at \url{https://github.com/SUSTech-Optimization/VF-iDCA}. The competitors are as follows:
\begin{itemize}\vspace{-5pt}
	\item \textbf{Implicit Differentiation:} We consider IGJO in \citet{feng2018gradient} and IFDM in \citet{bertrand2020implicit}. 
	\item \textbf{Grid Search:} We perform a brute-force grid search over a 10$\times$ 10 uniformly-spaced grid.
	\item \textbf{Random Search:} We apply uniform random sampling 100 times at each direction of hyperparameters.
	\item \textbf{TPE:} Tree-structured Parzen Estimator approach \cite{bergstra2013making} is a Bayesian optimization method based on Gaussian mixture models which can handle a large number of hyperparameters. 
\end{itemize}\vspace{-5pt}

All experiments run on a computer with Intel(R) Core(TM) i9-9900K CPU @ 3.60GHz and 16.00 GB memory.
We solve the strongly convex subproblem \eqref{p-DCA2_subproblem} at each iteration of VF-iDCA using the CVXPY package. 
IGJO and IFDM are implemented using code from \url{https://github.com/jjfeng/nonsmooth-joint-opt} and \url{https://github.com/QB3/sparse-ho}, respectively. 
We only use IFDM for elastic net, as its code can only deal with elastic net  among our tested problems.
Each subproblem in grid search and random search is solved using the CVXPY package.
TPE is implemented using code from \url{https://github.com/hyperopt/hyperopt} and its subproblem is solved using the CVXPY package. Whenever used, the CVXPY package is applied with the open source solvers ECOS and SCS only. 

As for the paramters $\delta_{\alpha}$ and $c_{\alpha}$ in VF-iDCA, we used $\delta_{\alpha} = 5$ for all the experiments, and we adopted $c_{\alpha} = 0.1$ in sparse group lasso while we used $c_{\alpha} = 1$ for the other applications. When $\delta_{\alpha}$ is too small, the algorithm may converge slower since the algorithm takes more steps to reach the exact penalty. However, in our experience, the performance of \cref{p-ipDCA} is not sensitive to the choice of $\delta_\alpha$. The choice of $c_\alpha$ is more important. We tested various ratios for $\Delta^k$ and $t^k$ to adjust $c_\alpha$. When $c_\alpha$ is large, the updating frequency of the penalty parameter will be high, which may slow down the convergence speed.

\subsection{Numerical experiments on synthetic data}

We consider hyperparameter selection for elastic net and variations on sparse group lasso, low-rank matrix completion, and support vector machines previously explored in \citet{feng2018gradient} and \citet{kunapuli2008classification}; see Section \ref{sectionC} of the Supplemental Materials for details. Detailed descriptions of the synthetic data generation settings are in Section \ref{sectionD} of the Supplemental Materials.

\subsubsection{Elastic Net}
\label{numerical_EN}

The numerical results on elastic net averaged over 30 data sets are in \cref{table:en_full_3}.  
Overall, VF-iDCA attains the highest-quality solutions relatively quickly. In the first two settings, the gradient-based methods (IGJO/IFDM and VF-iDCA) perform better than the gradient-free methods. Furthermore, our VF-iDCA achieves significantly lower validation and test errors than all other methods in the third setting, where $p = 2500$. This suggests that our method is especially competitive when the dimension of the problem is large. 

\begin{table}[th] 
\vspace{-5pt}
\caption{Elastic net problems on synthetic data. Here, $|I_{\text{tr}}|; |I_{\text{val}}|; |I_{\text{test}}|,$ and $p$ represent the number of training examples, the number of validation examples, the number of test examples, and the number of features, respectively. }\label{table:en_full_3}
\vspace{-5pt}
\begin{center}
\begin{small}
\resizebox{.49\textwidth}{!}{
\begin{tabular}{ll ccc}
\toprule
Settings & Method
& Time & Val. Err. & Test Err. \\
\midrule
\multirow{6}{*}{ \makecell[l]{ $|I_{\text{tr}}| = 100$\\$|I_{\text{val}}| = 20$\\
$|I_{\text{test}}| = 250$\\ $p = 250$ } } 
& Grid & 3.10 $\pm$ 0.44 & 6.16 $\pm$ 2.35 & 6.68 $\pm$ 1.16 \\
& Random & 3.55 $\pm$ 0.58 & 5.98 $\pm$ 2.24 & 6.67 $\pm$ 1.15 \\
& TPE & 5.41 $\pm$ 0.75 & 6.05 $\pm$ 2.30 & 6.77 $\pm$ 1.04 \\
& IGJO & 2.04 $\pm$ 1.46 & 4.43 $\pm$ 1.77 & 5.13 $\pm$ 1.37\\
& IFDM & 1.33 $\pm$ 0.55 & 4.41 $\pm$ 0.96 & 4.77 $\pm$ 1.46 \\
& VF-iDCA & \textbf{0.91 $\pm$ 0.19} & 1.95 $\pm$ 0.81 & \textbf{3.99 $\pm$ 0.69} \\
\midrule
\multirow{6}{*}{ \makecell[l]{ $|I_{\text{tr}}| = 100$\\$|I_{\text{val}}| = 100$\\
$|I_{\text{test}}| = 250$\\ $p = 250$ } }
& Grid & 3.17 $\pm$ 0.43 & 6.51 $\pm$ 1.53 & 6.82 $\pm$ 1.10 \\
& Random & 5.29 $\pm$ 0.60 & 6.44 $\pm$ 1.53 & 6.77 $\pm$ 1.14 \\
& TPE & 5.40 $\pm$ 0.84 & 6.44 $\pm$ 1.53 & 6.76 $\pm$ 1.06 \\
& IGJO & 2.42 $\pm$ 1.30 & 4.71 $\pm$ 1.32 & 4.88 $\pm$ 1.30 \\
& IFDM & 1.30 $\pm$ 0.41 & 4.78 $\pm$ 1.12 & 4.61 $\pm$ 1.12 \\
& VF-iDCA & \textbf{1.37 $\pm$ 0.29} & 3.04 $\pm$ 0.74 & \textbf{3.58 $\pm$ 0.60} \\
\midrule
\multirow{6}{*}{ \makecell[l]{ $|I_{\text{tr}}| = 100$\\$|I_{\text{val}}| = 100$\\
$|I_{\text{test}}| = 100$\\ $p = 2500$ } } 
& Grid & 19.05 $\pm$ 1.63 & 7.95 $\pm$ 1.10 & 8.54 $\pm$ 0.81 \\
& Random & 35.42 $\pm$ 3.55 & 7.90 $\pm$ 1.09 & 8.52 $\pm$ 0.79 \\
& TPE & 32.17 $\pm$ 7.40 & 7.89 $\pm$ 1.11 & 8.60 $\pm$ 0.87 \\
& IGJO & 16.12 $\pm$ 40.95 & 7.99 $\pm$ 1.18 & 8.41 $\pm$ 0.86 \\
& IFDM & \textbf{4.38 $\pm$ 2.53} & 7.97 $\pm$ 0.83 & 8.53 $\pm$ 1.53 \\
& VF-iDCA & 19.97 $\pm$ 5.17 & 1.61 $\pm$ 1.85 & \textbf{5.10 $\pm$ 1.07} \\
\bottomrule
\end{tabular}
}
\end{small}
\end{center}
\vspace{-15pt}
\end{table}

\subsubsection{Sparse Group Lasso}
\cref{table:sgl_full_3} records the numerical results of sparse group lasso averaged over 30 repetitions. Not only are the solutions of VF-iDCA are better than the gradient free methods, they are also better than the gradient-based IGJO. Furthermore, VF-iDCA is much faster than IGJO. 

\begin{table}[th]\vspace{-5pt}
\caption{Sparse group lasso problems on synthetic data. Recall that \#$\lambda$ is the dimension of hyperparameters, $p$ is the number of features and $M$ is the number of groups. Each group contains $p/M$ features. }
\label{table:sgl_full_3}
\vspace{-5pt}
\begin{center}
\begin{small}
\resizebox{.49\textwidth}{!}{
\begin{tabular}{ll rccc}
\toprule
Settings & Method
& \#$\lambda$ & Time & Val. Err. & Test Err. \\
\midrule
\multirow{5}{*}{\makecell[l]{$p=600$ \\ $M=30$}} 
& Grid & 2 & 30.38 $\pm$ 1.82 & 42.45 $\pm$ 7.67 & 44.56 $\pm$ 7.33 \\
& Random & 31 & 28.54 $\pm$ 1.51 & 39.27 $\pm$ 7.32 & 43.00 $\pm$ 8.83 \\
& TPE & 31 & 47.07 $\pm$ 4.01 & 35.69 $\pm$ 5.92 & 40.59 $\pm$ 6.67 \\
& IGJO & 31 & 69.62 $\pm$ 47.76 & 30.16 $\pm$ 7.41 & 39.28 $\pm$ 6.56 \\
& VF-iDCA & 31 & \textbf{8.13 $\pm$ 1.20} & 0.01 $\pm$ 0.00 & \textbf{38.50 $\pm$ 6.00} \\
\midrule
\multirow{5}{*}{\makecell[l]{$p=600$ \\ $M=300$}} 
& Grid & 2 & 20.84 $\pm$ 1.04 & 41.88 $\pm$ 7.64 & 44.90 $\pm$ 7.02 \\
& Random & 301 & 18.94 $\pm$ 1.09 & 43.92 $\pm$ 8.77 & 47.90 $\pm$ 8.55 \\
& TPE & 301 & 76.82 $\pm$ 2.55 & 39.22 $\pm$ 6.26 & 42.93 $\pm$ 8.00 \\
& IGJO & 301 & 160.85 $\pm$ 71.50 & 20.37 $\pm$ 4.46 & 38.52 $\pm$ 6.78 \\
& VF-iDCA & 301 & 56.73 $\pm$ 92.48 & 19.61 $\pm$ 8.33 & \textbf{33.55 $\pm$ 4.71} \\
\midrule
\multirow{5}{*}{\makecell[l]{$p=1200$ \\ $M=300$}} 
& Grid & 2 & 87.20 $\pm$ 5.85 & 49.56 $\pm$ 10.76 & 51.85 $\pm$ 12.90 \\
& Random & 301 & 73.75 $\pm$ 4.28 & 53.65 $\pm$ 12.03 & 55.84 $\pm$ 14.25 \\
& TPE & 301 & 117.07 $\pm$ 5.66 & 45.94 $\pm$ 9.30 & 51.67 $\pm$ 12.29 \\
& IGJO & 301 & 98.35 $\pm$ 47.47 & 20.70 $\pm$ 4.70 & 38.90 $\pm$ 7.20 \\
& VF-iDCA & 301 & \textbf{23.41 $\pm$ 1.31} & 17.90 $\pm$ 3.47 & \textbf{36.90 $\pm$ 7.48}\\
\bottomrule
\end{tabular}
}
\end{small}
\end{center}
\vspace{-20pt}
\end{table}

\subsubsection{Low-rank matrix completion}
For this problem, we present the numerical results on $60\times60$ matrices in \cref{table:MCG}. 

\begin{table}[hbt]
\vspace{-5pt}
\caption{Low-rank matrix completion problems on synthetic data. \#$\lambda$ is the dimension of hyperparameters.} \vspace{-10pt}
\label{table:MCG}
\begin{center}
\begin{small}
\resizebox{.49\textwidth}{!}{
\begin{tabular}{l rccc}
\toprule
Method
& \#$\lambda$ & Time & Val. Err. & Test Err. \\
\midrule
Grid & 2 & 20.67 $\pm$ 0.90 & 0.71 $\pm$ 0.21 & 0.76 $\pm$ 0.20 \\
Random & 25 & 32.49 $\pm$ 1.84 & 0.73 $\pm$ 0.21 & 0.80 $\pm$ 0.20 \\
TPE & 25 & 35.05 $\pm$ 9.37 & 0.68 $\pm$ 0.20 & 0.76 $\pm$ 0.18 \\
IGJO & 25 & 1268.65 $\pm$ 365.99 & 0.68 $\pm$ 0.21  &  0.72 $\pm$ 0.18  \\
VF-iDCA & 25 & 51.55 $\pm$ 10.43 & 0.06 $\pm$ 0.07 & \textbf{0.70 $\pm$ 0.16} \\
\bottomrule
\end{tabular}
}
\end{small}
\end{center}
\vspace{-20pt}
\end{table}

\subsection{Application to real data}

In this part, we will conduct numerical experiments on real datasets.
All datasets are from the LIBSVM repository \cite{libsvm}\footnote{\url{http://www.csie.ntu.edu.tw/~cjlin/libsvmtools/datasets/}.}. For each dataset, we will repeat the experiment 30 times. For each repetition of experiments, we applied a random shuffle to the original dataset before partition, and ran all the algorithms on the same partition.

\subsubsection{Elastic Net}
This section we will apply our method to choose hyperparameters for the elastic net in real data. We will test on three moderate scale datasets, including gisette \cite{guyon2004result}, duke breast-cancer \cite{west2001predicting} and sensit \cite{duarte2004vehicle}. 
For datasets gisette, duke breast-cancer, sensit, we randomly extracted 50, 11, 25 examples as training set, respectively; 50, 11, 25 examples as validation set, respectively; and the remaining for testing.

\begin{table}[th]
\vspace{-10pt}
\caption{Elastic net on madelon, gisette and duke breast-cancer. $p$, $|I_{\text{tr}}|$, $|I_{\text{val}}|$, $|I_{\text{test}}|$ denote features number, the size of train data, validation data and test data.}
\vspace{-10pt}
\label{table:realdata_EN_3}
\begin{center}
\begin{small}
\resizebox{.49\textwidth}{!}{
\begin{tabular}{ll ccc}
\toprule
Dataset & Method
& Time & Val. Err. & Test Err.\\
\midrule
\multirow{6}{*}{\makecell[l]{gisette \\$p=5000$\\ $|I_{\text{tr}}| = 50$\\$|I_{\text{val}}| = 50$\\ $|I_{\text{test}}| = 5900$}}
& Grid & 31.85 $\pm$ 3.51 & 0.22 $\pm$ 0.04 & 0.23 $\pm$ 0.01\\
& Random & 55.04 $\pm$ 9.05 & 0.22 $\pm$ 0.05 & 0.23 $\pm$ 0.02\\
& TPE & 39.12 $\pm$ 6.96 & 0.22 $\pm$ 0.05 & 0.24 $\pm$ 0.02\\
& IGJO & 6.10 $\pm$ 3.74 & 0.24 $\pm$ 0.05 & 0.25 $\pm$ 0.03\\
& IFDM & 191.18 $\pm$ 202.30 & 0.22 $\pm$ 0.02 & 0.23 $\pm$ 0.03 \\
& VF-iDCA & \textbf{5.39 $\pm$ 1.62} & 0.00 $\pm$ 0.00 & \textbf{0.19 $\pm$ 0.01}\\
\midrule
\multirow{6}{*}{\makecell[l]{duke breast-cancer \\$p=7129$\\ $|I_{\text{tr}}| = 11$\\$|I_{\text{val}}| = 11$\\ $|I_{\text{test}}| = 22$}}
& Grid & 23.75 $\pm$ 4.42 & 0.29 $\pm$ 0.06 & 0.36 $\pm$ 0.08\\
& Random & 39.68 $\pm$ 7.89 & 0.28 $\pm$ 0.07 & 0.41 $\pm$ 0.22\\
& TPE & 36.22 $\pm$ 22.16 & 0.27 $\pm$ 0.07 & 0.41 $\pm$ 0.24 \\
& IGJO & 3.51 $\pm$ 2.06 & 0.30 $\pm$ 0.06 & 0.37 $\pm$ 0.08\\
& IFDM & 26.79 $\pm$ 73.39 & 0.32 $\pm$ 0.05 & 0.36 $\pm$ 0.08 \\
& VF-iDCA & 12.23 $\pm$ 15.48 & 0.00 $\pm$ 0.00 & \textbf{0.29 $\pm$ 0.07}\\
\midrule 
\multirow{6}{*}{\makecell[l]{sensit \\$p=78823$\\ $|I_{\text{tr}}| = 25$\\$|I_{\text{val}}| = 25$\\ $|I_{\text{test}}| = 50$}}
& Grid & 1.26 $\pm$ 0.11 & 1.30 $\pm$ 0.82 & 1.31 $\pm$ 0.48 \\
& Random & 1.29 $\pm$ 0.11 & 1.52 $\pm$ 1.15 & 1.45 $\pm$ 0.60 \\
& TPE & 1.76 $\pm$ 0.09 & 1.38 $\pm$ 0.96 & 1.39 $\pm$ 0.55 \\
& IGJO & 0.49 $\pm$ 0.74 & 0.52 $\pm$ 0.20 & 0.61 $\pm$ 0.11 \\
& IFDM & 6.60 $\pm$ 3.49 & 0.57 $\pm$ 0.10 & 0.61 $\pm$ 0.23 \\
& VF-iDCA & \textbf{0.16 $\pm$ 0.08} & 0.23 $\pm$ 0.11 & \textbf{0.51 $\pm$ 0.06} \\
\bottomrule
\end{tabular}
}
\end{small}
\end{center}
\vspace{-10pt}
\end{table}

\subsubsection{Support Vector Machine}

We apply our experiments to moderate scale datasets, liver-disorders, diabetes, breast-cancer, sonar, a1a \cite{asuncion2007uci}, w1a \cite{Platt1999FastTO}. On each dataset, we choose hyperparameters for the SVM via $3/6$-fold cross-validation 30 times with the following data partition rule. For a dataset which contains $N$ examples, we extracted $3 \lfloor N/6 \rfloor$ examples as the training set used for cross-validation $\Omega$ and took the remaining part as the test set $\Omega_{\text{test}}$. 
Since the bilevel program for cross-validated SVM involves LL constraints, IGJO cannot be applied.
We implement VF-iDCA with two different tolerances. We denote the one with $tol=0.01$ by VF-iDCA, and the one with $tol=0.1$ by VF-iDCA-t.
We set $\bar{\mathbf{w}}_{lb} = 10^{-6}$ and $\bar{\mathbf{w}}_{ub} = 10$.

\begin{figure*}[htb]
	\centering
	\includegraphics[width=.95\textwidth, trim = 5cm 0cm 5cm 0cm, clip]{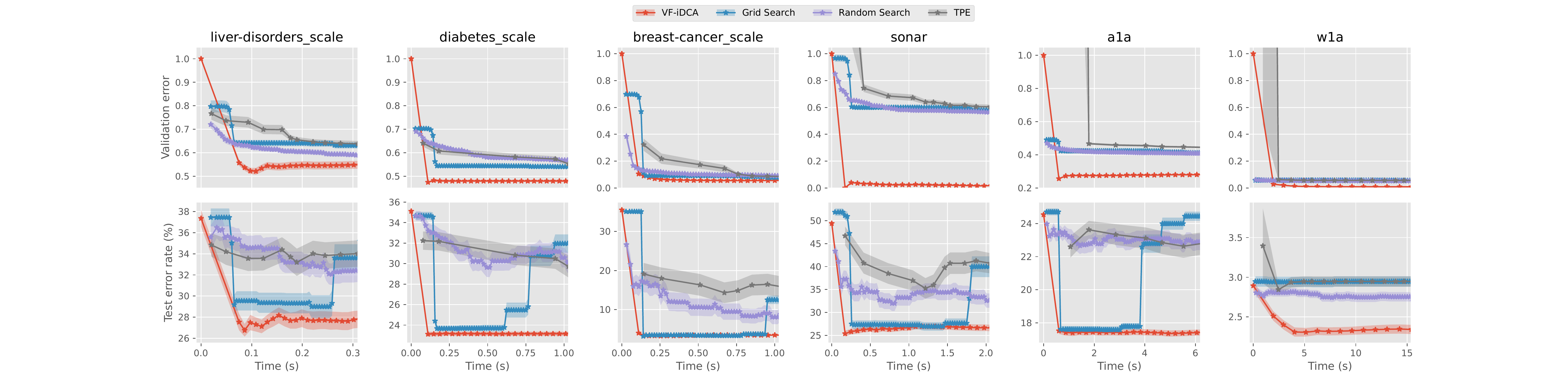}
	\vspace{-10pt}
	\caption{Comparison of the algorithms on SVM problem (validation error and test error versus time) for 6 datasets: liver-disorders\_scale, diabetes\_scale, breast-cancer\_scale, sonar, a1a, w1a}
	\label{fig:svm6}\vspace{-10pt} 
\end{figure*}

\cref{table:svm_full_3} and 
\cref{table:svm_full_6} record the numerical results of 3-fold SVM and 6-fold SVM on 6 datasets with 30 repeatations, respectively.
\cref{fig:svm6} shows how the validation error and test error rate changed over time for each algorithm. 
For such a problem, VF-iDCA shows superiority in both computation time and solution quality. By comparing the numerical results of VF-iDCA and VF-iDCA-t, it can be observed that VF-iDCA can obtain a satisfactory result with a moderate algorithmic tolerance.

The boundedness assumption for the sequence of adaptive penalty parameters, which plays an important role in the convergence analysis, appeared to be satisfied in our extensive experiments. In fact, as discussed in \citet{DCbilevel}, if we relax the constraint in \cref{p-singlelevel_problem} to be $l(x) - v(r) \leq \epsilon$ for any small fixed positive $\epsilon$, then the penalty parameters are guaranteed to be bounded. Thus, a small modification of Step 3 of Algorithm \ref{p-ipDCA} would solve this alternative ``$\epsilon$-problem”, and results analogous to those in Section \ref{sec:theory} would hold without assuming boundedness. In the future, we may also try to derive sufficient conditions that directly guarantee the boundedness of the adaptive penalty parameters for Algorithm \ref{p-ipDCA}.

\begin{table}[htbp]\vspace{-10pt}
\caption{3-fold SVM on datasets liver-disorders\_scale, diabetes\_scale , breast-cancer\_scale, sonar, a1a, w1a. \#$\lambda$ is the dimension of hyperparameters.}
\vspace{-10pt}
\label{table:svm_full_3}
	\begin{center}
\begin{small}
\resizebox{.49\textwidth}{!}{
\begin{threeparttable}
\begin{tabular}{ll rccc}
\toprule
Dataset & Method
& \#$\lambda$ & Time & Val. Err. & Test Err.\\
\midrule
\multirow{6}{*}{\makecell[l]{liver-disorders\_scale \\ $p=5$ \\ $|\Omega| = 72$ \\ $|\Omega_{\text{test}}| = 73$}}
& Grid & 2 & 0.53 $\pm$ 0.01 & 0.64 $\pm$ 0.08 & 0.34 $\pm$ 0.06 \\
& Random & 6 & 0.56 $\pm$ 0.03 & 0.58 $\pm$ 0.06 & 0.32 $\pm$ 0.05 \\
& TPE2 & 2 & 2.88 $\pm$ 1.16 & 0.61 $\pm$ 0.07 & 0.33 $\pm$ 0.06\\
& TPE & 6 
& 0.37 $\pm$ 0.29 & 0.65 $\pm$ 0.07 & 0.34 $\pm$ 0.05 \\
& VF-iDCA & 6 & 0.20 $\pm$ 0.06 & 0.53 $\pm$ 0.09 & 0.28 $\pm$ 0.05 \\
& VF-iDCA-t & 6 & \textbf{0.09 $\pm$ 0.02} & 0.53 $\pm$ 0.07 & \textbf{0.27 $\pm$ 0.03} \\
\midrule
\multirow{6}{*}{\makecell[l]{diabetes\_scale \\ $p=8$ \\ $|\Omega| = 384$ \\ $|\Omega_{\text{test}}| = 384$}}
& Grid & 2 & 1.70 $\pm$ 0.11 & 0.55 $\pm$ 0.03 & 0.33 $\pm$ 0.05 \\
& Random & 9 & 1.83 $\pm$ 0.09 & 0.56 $\pm$ 0.04 & 0.30 $\pm$ 0.06 \\
& TPE2 & 2 & 18.67 $\pm$ 7.84 & 0.54 $\pm$ 0.03 & 0.32 $\pm$ 0.06\\
& TPE & 9 
& 6.64 $\pm$ 4.30 & 0.55 $\pm$ 0.03 & 0.29 $\pm$ 0.05 \\
& VF-iDCA & 9 & 0.28 $\pm$ 0.03 & 0.48 $\pm$ 0.02 & \textbf{0.23 $\pm$ 0.01} \\
& VF-iDCA-t & 9 & \textbf{0.18 $\pm$ 0.02} & 0.48 $\pm$ 0.02 & \textbf{0.23 $\pm$ 0.01} \\
\midrule
\multirow{6}{*}{\makecell[l]{breast-cancer\_scale \\ $p=14$ \\ $|\Omega| = 338$ \\ $|\Omega_{\text{test}}| = 345$}}
& Grid & 2 & 1.63 $\pm$ 0.04 & 0.08 $\pm$ 0.01 & 0.12 $\pm$ 0.06 \\
& Random & 11 & 1.80 $\pm$ 0.03 & 0.09 $\pm$ 0.01 & 0.08 $\pm$ 0.09 \\
& TPE2 & 2 & 14.72 $\pm$ 6.02 & 0.07 $\pm$ 0.01 & 0.09 $\pm$ 0.10 \\
& TPE & 11 
& 9.14 $\pm$ 4.55 & 0.09 $\pm$ 0.01 & 0.10 $\pm$ 0.11 \\
& VF-iDCA & 11 & 1.12 $\pm$ 0.59 & 0.05 $\pm$ 0.01 & \textbf{0.03 $\pm$ 0.01} \\
& VF-iDCA-t & 11 & \textbf{0.14 $\pm$ 0.01} & 0.09 $\pm$ 0.01 & 0.04 $\pm$ 0.01 \\
\midrule
\multirow{6}{*}{\makecell[l]{sonar \\ $p=60$ \\ $|\Omega| = 102$ \\ $|\Omega_{\text{test}}| = 106$}}
& Grid & 2 & 3.19 $\pm$ 0.10 & 0.58 $\pm$ 0.08 & 0.40 $\pm$ 0.12 \\
& Random & 61 & 3.23 $\pm$ 0.06 & 0.54 $\pm$ 0.06 & 0.34 $\pm$ 0.10 \\
& TPE2 & 2 & 18.47 $\pm$ 6.84 & 0.57 $\pm$ 0.08 & 0.37 $\pm$ 0.13 \\
& TPE & 61 
& 40.77 $\pm$ 7.12 & 0.64 $\pm$ 0.10 & 0.41 $\pm$ 0.12 \\
& VF-iDCA & 61 & 2.22 $\pm$ 1.50 & 0.00 $\pm$ 0.00 & \textbf{0.24 $\pm$ 0.04} \\
& VF-iDCA-t & 61 & \textbf{0.48 $\pm$ 0.09} & 0.03 $\pm$ 0.02 & \textbf{0.24 $\pm$ 0.04} \\
\midrule
\multirow{6}{*}{\makecell[l]{a1a \\ $p=123$ \\ $|\Omega| = 801$ \\ $|\Omega_{\text{test}}| = 804$}}
& Grid & 2 & 8.04 $\pm$ 0.15 & 0.41 $\pm$ 0.02 & 0.24 $\pm$ 0.02 \\
& Random & 124 & 8.62 $\pm$ 0.30 & 0.41 $\pm$ 0.02 & 0.22 $\pm$ 0.03 \\
& TPE2 & 2 & 65.51 $\pm$ 16.24 & 0.41 $\pm$ 0.02 & 0.24 $\pm$ 0.02 \\
& TPE & 124 
& 176.59 $\pm$ 17.38 & 0.42 $\pm$ 0.03 & 0.23 $\pm$ 0.03 \\
& VF-iDCA & 124 & 10.17 $\pm$ 5.47 & 0.27 $\pm$ 0.02 & 0.18 $\pm$ 0.01 \\
& VF-iDCA-t & 124 & \textbf{1.10 $\pm$ 0.07} & 0.27 $\pm$ 0.02 & \textbf{0.17 $\pm$ 0.01} \\
\midrule
\multirow{6}{*}{\makecell[l]{w1a \\ $p=300$ \\ $|\Omega| = 1236$ \\ $|\Omega_{\text{test}}| = 1241$}}
& Grid & 2 & 20.21 $\pm$ 0.82 & 0.06 $\pm$ 0.01 & 0.03 $\pm$ 0.00 \\
& Random & 301 & 20.44 $\pm$ 1.10 & 0.06 $\pm$ 0.01 & 0.03 $\pm$ 0.00 \\
& TPE2 & 2 & 86.10 $\pm$ 28.19 & 0.06 $\pm$ 0.01 & 0.03 $\pm$ 0.00 \\
& TPE & 301 
& 299.62 $\pm$ 78.72 & 0.06 $\pm$ 0.01 & 0.03 $\pm$ 0.00 \\
& VF-iDCA & 301 & 27.49 $\pm$ 7.31 & 0.01 $\pm$ 0.00 & \textbf{0.02 $\pm$ 0.00} \\
& VF-iDCA-t & 301 & \textbf{4.87 $\pm$ 0.51} & 0.01 $\pm$ 0.00 & \textbf{0.02 $\pm$ 0.00} \\
\bottomrule
\end{tabular}
\end{threeparttable}}
\end{small}
\end{center}
\vspace{-10pt}
\end{table}

\begin{table}[htbp]\vspace{-15pt}
	\caption{6-fold SVM on datasets liver-disorders\_scale, diabetes\_scale , breast-cancer\_scale, sonar, a1a, w1a. \#$\lambda$ is the dimension of hyperparameters.}
	\label{table:svm_full_6}
	\vspace{-5pt}
	\begin{center}
	\begin{small}
	\resizebox{.49\textwidth}{!}{
	\begin{threeparttable}
	\begin{tabular}{ll rccc}
		\toprule
		Dataset & Method
		& \#$\lambda$ & Time & Val. Err. & Test Err.\\
		\midrule
		\multirow{6}{*}{\makecell[l]{liver-disorders\_scale \\ $p=5$ \\ $|\Omega| = 72$ \\ $|\Omega_{\text{test}}| = 73$}}
		& Grid & 2 
		& 0.78 $\pm$ 0.02 & 0.63 $\pm$ 0.08 & 0.33 $\pm$ 0.07 \\
		& Random & 6 
		& 0.79 $\pm$ 0.04 & 0.62 $\pm$ 0.07 & 0.31 $\pm$ 0.05 \\
		& TPE2 & 2 
		& 6.88 $\pm$ 4.15 & 0.62 $\pm$ 0.07 & 0.32 $\pm$ 0.06 \\
		& TPE & 6 
		& 1.06 $\pm$ 1.04 & 0.63 $\pm$ 0.08 & 0.34 $\pm$ 0.05 \\
		& VF-iDCA & 6 
		& 0.33 $\pm$ 0.14 & 0.41 $\pm$ 0.08 & \textbf{0.27 $\pm$ 0.05} \\
		& VF-iDCA-t & 6
		& \textbf{0.18 $\pm$ 0.03} & 0.41 $\pm$ 0.08 & \textbf{0.27 $\pm$ 0.04} \\
		\midrule
		\multirow{6}{*}{\makecell[l]{diabetes\_scale \\ $p=8$ \\ $|\Omega| = 384$ \\ $|\Omega_{\text{test}}| = 384$}}
		& Grid & 2 
		& 3.18 $\pm$ 0.14 & 0.55 $\pm$ 0.03 & 0.32 $\pm$ 0.05 \\
		& Random & 9 
		& 3.63 $\pm$ 0.21 & 0.56 $\pm$ 0.03 & 0.31 $\pm$ 0.05 \\
		& TPE2 & 2 
		& 51.85 $\pm$ 20.49 & 0.55 $\pm$ 0.03 & 0.33 $\pm$ 0.05 \\
		& TPE & 9 
		& 29.52 $\pm$ 13.13 & 0.55 $\pm$ 0.03 & 0.27 $\pm$ 0.06 \\
		& VF-iDCA & 9 
		& 0.56 $\pm$ 0.08 & 0.43 $\pm$ 0.02 & \textbf{0.23 $\pm$ 0.01} \\
		& VF-iDCA-t & 9 
		& \textbf{0.35 $\pm$ 0.01} & 0.43 $\pm$ 0.02 & \textbf{0.23 $\pm$ 0.01} \\
		\midrule
		\multirow{6}{*}{\makecell[l]{breast-cancer\_scale \\ $p=14$ \\ $|\Omega| = 336$ \\ $|\Omega_{\text{test}}| = 347$}}
		& Grid & 2 
		& 3.38 $\pm$ 0.25 & 0.08 $\pm$ 0.02 & 0.15 $\pm$ 0.06 \\
		& Random & 11 
		& 3.92 $\pm$ 0.29 & 0.09 $\pm$ 0.02 & 0.07 $\pm$ 0.08 \\
		& TPE2 & 2 
		& 38.69 $\pm$ 16.27 & 0.07 $\pm$ 0.02 & 0.08 $\pm$ 0.09 \\
		& TPE & 11 
		& 25.96 $\pm$ 12.95 & 0.09 $\pm$ 0.01 & 0.11 $\pm$ 0.13 \\
		& VF-iDCA & 11 
		& 2.01 $\pm$ 0.17 & 0.03 $\pm$ 0.01 & \textbf{0.03 $\pm$ 0.01} \\
		& VF-iDCA-t & 11 
		& \textbf{0.29 $\pm$ 0.08} & 0.08 $\pm$ 0.01 & \textbf{0.04 $\pm$ 0.01}\\
		\midrule
		\multirow{6}{*}{\makecell[l]{sonar \\ $p=60$ \\ $|\Omega| = 102$ \\ $|\Omega_{\text{test}}| = 106$}}
		& Grid & 2 
		& 6.57 $\pm$ 0.32 & 0.59 $\pm$ 0.08 & 0.39 $\pm$ 0.11 \\
		& Random & 61 
		& 6.44 $\pm$ 0.28 & 0.54 $\pm$ 0.06 & 0.32 $\pm$ 0.08 \\
		& TPE2 & 2 
		& 58.19 $\pm$ 29.60 & 0.57 $\pm$ 0.08 & 0.36 $\pm$ 0.12 \\
		& TPE & 61 
		& 97.65 $\pm$ 31.37 & 0.60 $\pm$ 0.07 & 0.39 $\pm$ 0.12 \\
		& VF-iDCA & 61 
		& 0.92 $\pm$ 0.02 & 0.00 $\pm$ 0.00 & \textbf{0.23 $\pm$ 0.04} \\
		& VF-iDCA-t & 61 
		& \textbf{0.92 $\pm$ 0.02} & 0.00 $\pm$ 0.00 & \textbf{0.23 $\pm$ 0.04} \\
		\midrule
		\multirow{6}{*}{\makecell[l]{a1a \\ $p=123$ \\ $|\Omega| = 798$ \\ $|\Omega_{\text{test}}| = 807$}}
		& Grid & 2 
		& 17.60 $\pm$ 0.36 & 0.40 $\pm$ 0.02 & 0.25 $\pm$ 0.02 \\
		& Random & 124 
		& 18.59 $\pm$ 0.42 & 0.40 $\pm$ 0.02 & 0.21 $\pm$ 0.03 \\
		& TPE2 & 2 
		& 161.68 $\pm$ 42.67 & 0.40 $\pm$ 0.02 & 0.24 $\pm$ 0.02 \\
		& TPE & 124 
		& 312.63 $\pm$ 60.60 & 0.41 $\pm$ 0.03 & 0.23 $\pm$ 0.03 \\
		& VF-iDCA & 124 
		& 63.01 $\pm$ 186.14 & 0.18 $\pm$ 0.02 & \textbf{0.17 $\pm$ 0.01} \\
		& VF-iDCA-t & 124 & \textbf{4.22 $\pm$ 0.37} & 0.19 $\pm$ 0.01 & \textbf{0.17 $\pm$ 0.01} \\
		\midrule
		\multirow{6}{*}{\makecell[l]{w1a \\ $p=300$ \\ $|\Omega| = 1236$ \\ $|\Omega_{\text{test}}| = 1241$}}
		& Grid & 2 
		& 44.29 $\pm$ 1.39 & 0.05 $\pm$ 0.00 & 0.03 $\pm$ 0.00 \\
		& Random & 301 
		& 61.80 $\pm$ 2.91 & 0.05 $\pm$ 0.00 & 0.03 $\pm$ 0.00 \\
		& TPE2 & 2 
		& 190.04 $\pm$ 39.00 & 0.05 $\pm$ 0.00 & 0.03 $\pm$ 0.00 \\
		& TPE & 301 
		& 703.72 $\pm$ 82.75 & 0.05 $\pm$ 0.01 & 0.03 $\pm$ 0.00 \\
		& VF-iDCA & 301 
		& 97.50 $\pm$ 35.99 & 0.01 $\pm$ 0.00 & \textbf{0.02 $\pm$ 0.00} \\
		& VF-iDCA-t & 301 & \textbf{26.74 $\pm$ 3.67} & 0.01 $\pm$ 0.00 & \textbf{0.02 $\pm$ 0.00} \\
		\bottomrule
	\end{tabular}
	\end{threeparttable}}
	\end{small}
	\end{center}
	\vspace{-20pt}
\end{table}

\section{Conclusion}

\begin{figure}[htbp]\vspace{-20pt} 
	\centering 
	\includegraphics[width = .49\textwidth, trim = 2cm 0cm 1cm 0cm, clip]{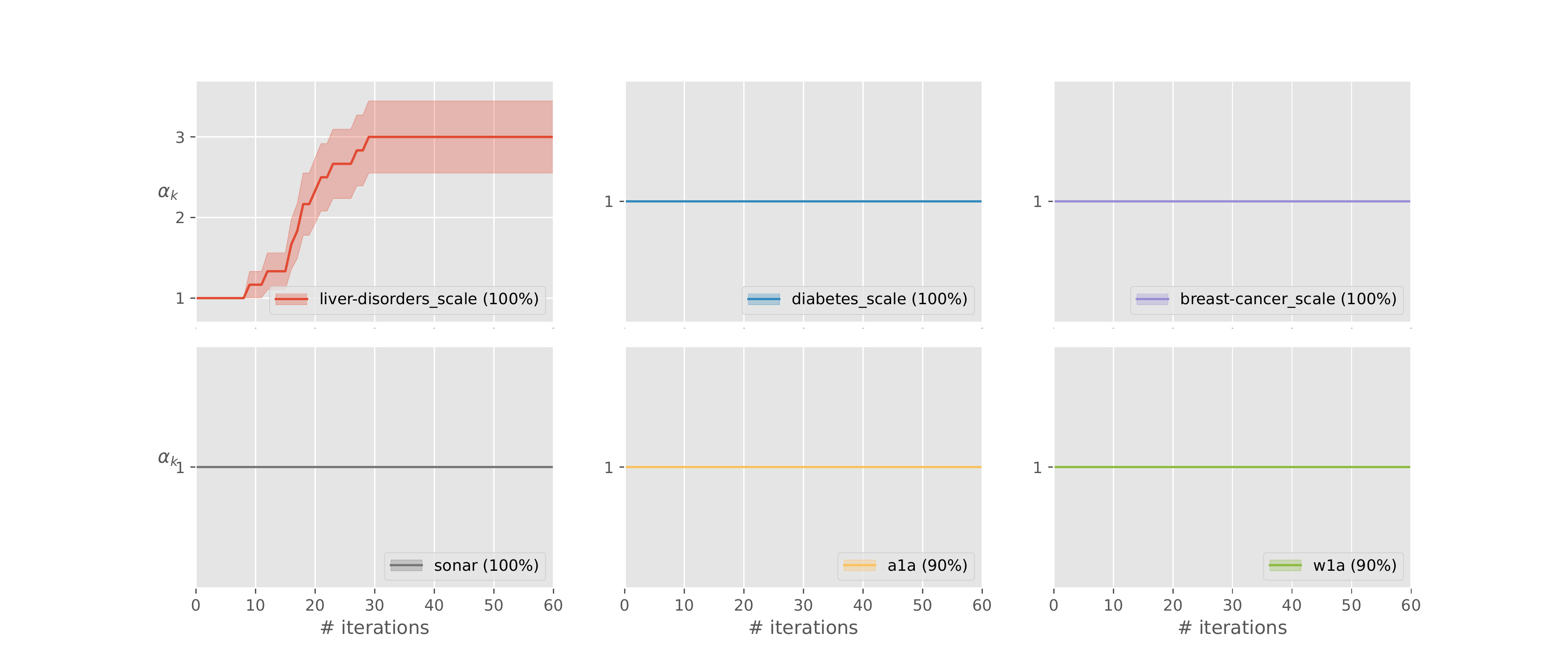}
	\vspace{-20pt}
	\caption{Adaptive penalty parameter $\alpha_k$ (defined in \cref{adaptpen})  of the VF-iDCA on SVM problem for 6 datasets: liver-disorders\_scale, diabetes\_scale, breast-cancer\_scale, sonar, a1a, w1a}
	\label{fig:beta}
	\vspace{-15pt}
\end{figure}

In this paper, we develop VF-iDCA, a new algorithm for hyperparameter tuning that can be applied for a single training/validation split or for cross-validation. This algorithm heavily exploits the structure of the hyperparameter tuning BLP in order to develop a tailored and effective solution.

Critically, unlike existing gradient-based BLP algorithms, we are able to rigorously establish  the first result for sequential convergence towards stationary solutions in the absence of LLSC and LLS simplifications. Furthermore, in our numerical experiments we find that VF-iDCA is able to match or outperform the state-of-the-art hyperparameter tuning approaches in terms of validation and test error for a number of machine learning algorithms, especially when the number of hyperparameters is large. In many cases, VF-iDCA also enjoys considerable gains in computational efficiency compared to the state-of-the-art.  

The boundedness assumption of the sequence of adaptive penalty parameters, which plays an important role in the convergence analysis, appeared to be satisfied in our extensive experiments. In fact, as discussed in \citet{DCbilevel}, if we relax the constraint in \cref{p-singlelevel_problem} to be $l(x) - v(r) \leq \epsilon$ for any small fixed positive $\epsilon$, then the penalty parameters are guaranteed to be bounded. Thus, a small modification of Step 3 of Algorithm \ref{p-ipDCA} would solve this alternative ``$\epsilon$-problem”, and results analogous to those in Section \ref{sec:theory} would hold without assuming boundedness. In the future, we may also try to derive sufficient conditions that directly guarantee the boundedness of the adaptive penalty parameters for Algorithm \ref{p-ipDCA}.

\vspace{-5pt}
\section*{Acknowledgements}
\vspace{-5pt}
The alphabetical order of the authors indicates the equal contribution to the paper. This work is partially supported by NSERC, the National Natural Science Foundation of China (No. 11971220), the Shenzhen Science and Technology Program (No. RCYX20200714114700072), the Guangdong Basic and Applied Basic Research Foundation (No. 2022B1515020082) and the Pacific Institute for the Mathematical Sciences (PIMS).


\clearpage
\bibliography{reference}
\bibliographystyle{icml2022}

\clearpage
\appendix
The supplemental materials are organized as follows. In \cref{sectionA}, we present all the detailed proofs of the theoretical results given in Section 3. And in \cref{sectionB}, the detailed description of the VF-iDCA for bilevel problem in general setting is presented. \cref{sectionC} is devoted to the illustration of hyperparameter decoupling for BLP applications.
Additional information on how the numerical experiments were run are given in \cref{sectionD}.

\section{Detailed Proofs}\label{sectionA}
\subsection{Proof of \cref{p-inclusion1}}
By  convexity, $x$ is an optimal solution for problem \cref{convexop1}  if and only if its first-order optimality condition  
\begin{equation}
  0\in \partial l(x)+\sum_{i=1}^J \lambda_i \partial P_i(x) \label{KKT1}
\end{equation} holds.
Assuming the  KKT condition for problem \cref{convexop2}
\begin{equation}\label{KKT2}
	\begin{aligned}
		\exists \lambda \in \mathbb{R}^J_+ \mbox{ s.t. }  & 0\in \partial l(x)+\sum_{i=1}^J \lambda_i \partial P_i(x), \\& \lambda_i (P_i(x)-r_i)=0, i=1,\dots, J, \\
		& P_i(x)-r_i \le 0, i=1,\dots, J,
	\end{aligned}
\end{equation}
holds, then $x\in S_c(r)$.
Since the difference between the KKT conditions \cref{KKT1} and \cref{KKT2} is the complementary slackness conditions $ \lambda_i (P_i(x)-r_i)=0, P_i(x)-r_i \le 0, i=1,\dots, J,$ which holds automatically if $P(x)=r$, the conclusion follows immediately.

\subsection{Proof of \cref{p-inclusion2}}
Since $\lambda \in \mathcal{M}(x,r)$, we have 
\[
0 \in \partial l(x) + \sum_{i = 1}^{J}\lambda_i \partial P_i(x).
\]
By convexity we have
$ x \in  S_{p}(\lambda)$.

\subsection{Proof of \cref{p-globalequi}}
Problems \cref{original_problem} and \cref{reform_problem} can be equivalently rewritten in the  following form
\begin{equation}
	\begin{aligned}\label{equiv(1)}
		\min ~~& L(x) \\
		\text{s.t.} ~~& x\in S_{p}(\lambda), \lambda \in \mathbb{R}_+^J,
	\end{aligned}
\end{equation}
and 
\begin{equation}
	\begin{aligned}\label{equiv(2)}
		\min ~~& L(x) \\
		\text{s.t.} ~~& x\in  S_c(r), r \in \mathbb{R}_+^J ,
	\end{aligned}
\end{equation}
respectively. 

(1) Suppose that $(\bar{x}, \bar{r}) \in \mathbb{R}^n \times \mathbb{R}_+^J$  is a global optimal solution of problem \cref{equiv(2)} and $\bar \lambda \in \mathcal{M}(\bar x,\bar r)$.
We want to show that $(\bar{x}, \bar \lambda)$ is a global optimal solution of problem \cref{equiv(1)}. Let $(x, \lambda)$ be a feasible solution of problem \cref{equiv(1)} which means that $x\in S_{p}(\lambda), \lambda \in \mathbb{R}_+^J$. By \cref{p-inclusion1}, $x\in S_c(r)$ with $r=P(x)\in  \mathbb{R}_+^J$ and hence $(x, r)$ is a feasible solution of problem \cref{equiv(2)}. Hence by the optimality of $(\bar{x}, \bar{r})$ to problem \cref{equiv(2)}, we have
$L(\bar x) \leq L(x)$. This shows that $(\bar{x}, \bar \lambda)$ is a global optimal solution of problem \cref{equiv(1)}.

(2) Suppose that $(\bar{x}, \bar{\lambda}) \in \mathbb{R}^n \times \mathbb{R}_+^J$  is a global optimal solution of problem \cref{equiv(1)}.
We want to show that $(\bar{x}, \bar r)$ with $\bar r:=P(\bar{x})$ is a global optimal solution of problem \cref{equiv(2)}. Since $\bar x\in S_p(\bar \lambda)$, by \cref{p-inclusion1}, $\bar x\in S_c(\bar r)$ with $\bar r=P(\bar x)$. Hence $(\bar x, \bar r)$ is a feasible solution of problem \cref{equiv(2)}. Let $(x,  r)$ be a feasible solution of problem \cref{equiv(2)} which means that $x\in S_{c}(r), r \in \mathbb{R}_+^J$. 
As $D$ is a dense subset of $X \times \mathbb{R}_+^J$, there exists a sequence $\{(x^k, r^k)\} \subseteq D$ satisfying $(x^k, r^k) \rightarrow (x, r)$.
Since $\mathcal{M}(x^k, r^k)\not =\emptyset$ for all $(x^k, r^k) \in D$, we can pick  $\lambda^k \in \mathcal{M}(x^k,r^k)$.
By \cref{p-inclusion2}, $x^k\in S_p(\lambda^k)$ and hence $(x^k, \lambda^k)$ is a feasible solution of problem \cref{equiv(1)}. By the optimality of $(\bar{x}, \bar{\lambda})$ to problem \cref{equiv(1)}, we have
$L(\bar x) \leq L(x^k)$. By taking $k \rightarrow \infty$ and the continuity of $L$, this shows that $(\bar{x}, \bar r)$ is a global optimal solution of problem \cref{equiv(2)}.

\subsection{Proof of \cref{local_min1}}
Since  $(\bar{x}, \bar r)$ is a local optimal solution of problem \cref{equiv(2)}, there is  a $\epsilon_0>0$ such that $\forall x\in  S_c(r), r \in \mathbb{R}_+^J$, and $(x, r) \in \mathbb{B}_{\epsilon_0}(\bar x, \bar{r})$,
\begin{equation} \label{localop1} L(\bar x) \leq L(x),
\end{equation} where $\mathbb{B}_{\epsilon_0}(\bar x,\bar r)$ denotes the closed  ball centered at $(\bar x, \bar{r})$ with radius $\epsilon_0$.
Let $(x, \lambda)$ be a feasible solution of problem \cref{equiv(1)}  which means that $x\in S_{p}(\lambda), \lambda \in \mathbb{R}_+^J$.  
By \cref{p-inclusion1}, $x\in S_c(r)$ with $r=P(x)\in  \mathbb{R}_+^J$ and hence $(x, r)$ is a feasible solution of problem \cref{equiv(2)}. Moreover since $P$ is continuous, we can find  $0<\epsilon_1<\epsilon_0$ such that as $x \in \mathbb{B}_{\epsilon_1}(\bar x)$, $r = P(x) \in \mathbb{B}_{\epsilon_0}(\bar{r})$.  Hence by \cref{localop1}, the local optimality of $(\bar{x}, \bar{r})$ to problem \cref{equiv(2)}, we have
$L(\bar x) \leq L(x)$. This shows that $(\bar{x}, \bar \lambda)$ is a local optimal solution of problem \cref{equiv(1)}. 

\subsection{Proof of \cref{local_min2}}

Since  $(\bar{x}, \bar \lambda)$ is a local optimal solution of problem \cref{equiv(1)} locally with respect to variable $x$ and globally with respect to variable $\lambda$, there is a $\epsilon_0>0$ such that  $\forall x\in  S_p(\lambda), \lambda \in \mathbb{R}_+^J$, and $x \in \mathbb{B}_{\epsilon_0}(\bar x)$,
\begin{equation} \label{localop2} L(\bar x) \leq L(x).
\end{equation} Since $\bar x\in S_p(\bar \lambda)$, by \cref{p-inclusion1}, $\bar x\in S_c(\bar r)$ with $\bar r=P(\bar x)$. Hence $(\bar x,\bar r)$ is a feasible solution of problem \cref{equiv(2)}. Let $(x, r)$ be a feasible solution of problem \cref{equiv(2)}) which means that $x\in S_{c}(r), r \in \mathbb{R}_+^J$.  
As $D$ is a dense subset of $X\times \mathbb{R}_+^J$, there exists a sequence $\{(x^k,r^k)\} \subseteq D$ satisfying $(x^k,r^k) \rightarrow (x,r)$.
Since $\mathcal{M}(x^k,r^k)\not =\emptyset$ for all $(x^k,r^k) \in D$, we can pick  $\lambda^k \in \mathcal{M}(x^k,r^k)$.
By \cref{p-inclusion2}, $x^k\in S_p(\lambda^k)$ and hence $(x^k, \lambda^k)$ is a feasible solution of problem \cref{equiv(1)}. 
Moreover suppose that $x$ and $\{x^k\}$ lie in $\mathbb{B}_{\epsilon_0}(\bar x,\bar u)$. Hence by (\ref{localop2}), we have
$L(\bar x) \leq L(x^k)$. Taking $k \rightarrow \infty$ implies $L(\bar x) \leq L(x)$. This shows that $(\bar{x}, \bar r)$ is a local optimal solution of problem \cref{equiv(2)} locally with respect to variable $x$ and globally with respect to variable $r$. 

Now suppose that $\mathcal{M}(\bar x,\bar r)=\{\bar \lambda\}$ is a singleton. Let $(\bar{x}, \bar \lambda)$ be a local optimal solution of problem \cref{equiv(1)}. Then there is  a $\epsilon_1>0$ such that $\forall x\in  S_p(\lambda), \lambda \in \mathbb{R}_+^J $ and $ (x, \lambda) \in \mathbb{B}_{\epsilon_1}(\bar x, \bar \lambda) $,
\begin{equation} \label{localop3} L(\bar x) \leq L(x).
\end{equation} 
Let $(x, r)$ be a feasible solution of problem \cref{equiv(2)} which lies in $\mathbb{B}_{\epsilon_2}(\bar x, \bar r)$ for some $0<\epsilon_2<\epsilon_1$. It implies that $x\in S_{c}(r),  r \in \mathbb{R}_+^J$ and $(x,r)\in \mathbb{B}_{\epsilon_2}(\bar x,\bar r)$.  
By the assumption of the uniqueness of the multipliers, we can show through proof by contradiction that when $\epsilon_2$ is sufficiently small, for $x \in \mathbb{B}_{\epsilon_2}(\bar x)$, $\Lambda \cap \mathcal{M}(x,r) \subseteq \mathbb{B}_{\epsilon_1}(\bar \lambda)$.	
We can find a sequence $\{(x^k,r^k)\} \subseteq D \cap \mathbb{B}_{\epsilon}(\bar x, \bar{r})$ satisfying $(x^k, r^k) \rightarrow (x, r)$.
Since $\Lambda \cap \mathcal{M}(x^k, r^k)\not =\emptyset$ for all $(x^k, r^k) \in D \cap \mathbb{B}_{\epsilon}(\bar x, \bar{r})$, we can pick  $\lambda^k \in \Lambda \cap \mathcal{M}(x^k,r^k)$. We can assume without loss of generality that $(x^k) \in \mathbb{B}_{\epsilon_2}(\bar x)$ and thus $\lambda^k \in \mathbb{B}_{\epsilon_1}(\bar \lambda)$.
By \cref{p-inclusion2}, $x^k\in S_p(\lambda^k)$ and hence $(x^k, \lambda^k)$ is a feasible solution of problem \cref{equiv(1)}. 
Hence by \cref{localop3}, we have $L(\bar x) \leq L(x_k)$. Taking $k \rightarrow \infty$ gives us $L(\bar x) \leq L(x)$. This shows that $(\bar{x}, \bar r)$ is a local optimal solution of problem \cref{equiv(2)}.

\subsection{Proof of \cref{p-suff_decreasenew}}

For convenience,  in the proof we denote by $g_k(z):=\displaystyle \max_{ i=1,\ldots,J}\left \{0, V_k(z), P_i(x) - r_i\right \}$. Then
\begin{equation}
\phi_k(z) = L(x)+\frac{\rho}{2}\|z-z^k\|^2 +\alpha_k g_k(z),\label{phi} \end{equation}
where $V_k(z)=l(x)-v(r^k)+\langle \gamma^k, r-r^k\rangle$.   By using the same arguments as in the proof of Lemma 1 in \cite{DCbilevel}, we can have following decrease result on $\phi_k$,
\begin{equation}\label{suff_decrease_proof_eq1}
	\begin{aligned}
		\phi_k(z^{k+1}) \le \phi_k(z^{k}) + \frac{\rho}{4}\|z^{k} - z^{k-1}\|^2 .
	\end{aligned}
\end{equation}
Next since  $v^*$ is the conjugate function of $v$,  by definition
\begin{equation}
	-v(r^k) \le  \langle \gamma^{k-1}, r^k \rangle + v^*( - \gamma^{k-1}).\label{conjug}
\end{equation}
Denote by $$\tilde g(z,\gamma):=\max_{ i=1,\ldots,J}\{0, l(x)+\langle \gamma, r\rangle +v^*(-\gamma), P_i(x) - r_i\}.$$ Then
$$E_{\alpha} (z,z_0,\gamma):=L(x)+\frac{\rho}{4}\|z-z_0\|^2 +\delta_\Sigma (z)+\alpha_k \tilde g(z,\gamma).$$

Combining \cref{phi}, \cref{suff_decrease_proof_eq1} and \cref{conjug},  and taking into account the fact that $z^k\in \Sigma$,  we obtain
\begin{equation}\label{suff_decrease_proof_eq2}
	{\phi}_k(z^{k+1}) \le 	E_{\alpha_k}(z^k,z^{k-1},\gamma^{k-1}). 
\end{equation}
Again, as $v^*$ is the conjugate function of $v$, and $- \gamma^k \in \partial v(r^k)$, there holds that
\begin{equation*}
	-v(r^k) - \langle \gamma^{k}, r^k \rangle = v^*(-\gamma^{k}),
\end{equation*}
which implies
\begin{equation}\label{suff_decrease_proof_eq3}
	{\phi}_k(z^{k+1}) = E_{\alpha_k}(z^{k+1},z^{k},\gamma^{k}) + \frac{\rho}{4} \|z^{k+1} - z^k\|^2.
\end{equation}
Combining with \cref{suff_decrease_proof_eq1} and \cref{suff_decrease_proof_eq2} gives us \cref{p-suff_decrease_eq2}.

It remains to prove \cref{p-re_err_eq}.
Since  $l(x) $ and $P_i(x) - {r}_i$, $i = 1, \ldots, J$ are all convex and continuous  and thus regular,  by the calculus rule for the pointwise maximum (see e.g.  Proposition 2.3.12 in \cite{clarke1990optimization}),  we have that $g_k(z)$ is regular and
\begin{equation}
	\begin{aligned}
			&\partial g_k(z) = \\
		&  \left  \{ \eta ( \partial l(x) \times \{\gamma^k\}) +\sum_{i=1}^J \lambda_i \partial P_i(x)\times \{-\lambda_i e_i\} \right .\\
		& \left .
		\mbox{ s.t. }  \lambda_i \in [0,1],  \lambda_i \left (P_i(x)-r_i-g_k(z) \right )=0, \right . \\
		& \left .  \eta \in [0,1], \eta (l(x)-v(r^k)+\langle \gamma^k, r-r^k\rangle-g_k(z))=0,\right.\\
		& \left .g_k(z)(1-\eta-\sum_{i=1}^J \lambda_i) =0, \eta+\sum_{i=1}^J \lambda_i\leq 1 \right \},
	\end{aligned}\label{g1}
	\end{equation}
	where $e_i$ denotes the unit vector with the $i$th component  equal to $1$.

Similarly,
$\tilde g(z,\gamma)$ is regular and 
\begin{equation}
	\begin{aligned}
		& \partial_z \tilde g(z,\gamma^k) =  \\
		&  \left  \{ \eta (\partial l(x) \times \{\gamma^k\}) +\sum_{i=1}^J \lambda_i \partial P_i(x)\times \{-\lambda_i e_i\} \right .\\
		& \left .
		\mbox{ s.t. }  \lambda_i \in [0,1],  \lambda_i \left (P_i(x)-r_i-\tilde g(z,\gamma^k) \right )=0, \right . \\
		& \left .  \eta \in [0,1], \eta (l(x)+\langle \gamma^k, r\rangle +v^*(-\gamma^k)-\tilde g(z,\gamma^k))=0, \right.\\
		& \left .\tilde g(z,\gamma^k)(1-\eta-\sum_{i=1}^J \lambda_i) =0, \eta+\sum_{i=1}^J \lambda_i\leq 1 \right \}.
	\end{aligned}\label{g2}
\end{equation}
We also have 
\begin{equation}
	\begin{aligned}
		& \partial_\gamma  \tilde g(z^{k+1},\gamma) =   \left  \{ \eta (r^{k+1}-\partial v^*(-\gamma)) \mbox{ s.t. }   \eta \in [0,1], \right .\\
		& \left. \eta \left  (l(x)+\langle \gamma, r^{k+1}\rangle +v^*(-\gamma)-\tilde g(z^{k+1},\gamma)\right )=0 \right \}.
	\end{aligned}\label{g3}
\end{equation}

Since $v^*$ is the conjugate function of $v$, and $- \gamma^k \in \partial v(r^k)$, there holds that
\begin{equation*}
	-v(r^k) - \langle \gamma^{k}, r^k \rangle = v^*(-\gamma^{k}).
\end{equation*}
Hence 
\[
\begin{aligned}
V_k(z^{k+1}) = \, &l(x^{k+1})-v(r^k)+\langle \gamma^k, r^{k+1}-r^k\rangle \\
 =\, & l(x^{k+1}) + \langle \gamma^k, r^{k+1}\rangle + v^*(- \gamma^k), 
\end{aligned}
\]
and thus
\[
 g_k(z^{k+1}) = \tilde g(z^{k+1},\gamma^k),
\]
and 
\begin{equation}\label{subdiff_equal}
	\partial g_k(z^{k+1}) = \partial_z \tilde  g(z^{k+1},\gamma^k).
\end{equation}
Since  $L(x) $ and $g_k(z)$ are all convex and continuous  and thus regular,   according to the the subdifferential sum rule we have
 (see e.g. Proposition 2.3.3  in \citet{clarke1990optimization}),  we have
$$\partial \phi_k(z) =\partial L(x)\times \{0\} +{\rho} (z-z^k) + \alpha_k \partial g_k(z).$$
Similarly we have
$$
\begin{aligned}
\partial_z E_{\alpha_k} (z, z^k,\gamma^k) = \,&\partial L(x)\times \{0\} +\frac{\rho}{2} (z-z^k) \\
&+ \mathcal{N}_\Sigma(z^{k+1})+\alpha_k \partial_z \tilde  g(z,\gamma^k).
\end{aligned}
$$
Hence by \cref{subdiff_equal},  we have
\begin{equation}
\begin{aligned}
	&\partial_z E_{\alpha_k}(z^{k+1},z^{k},\gamma^{k}) \\
	=\,&\partial \phi_k(z^{k+1} )+ \mathcal{N}_\Sigma(z^{k+1}) - \frac{\rho}{2}(z^{k+1}-z^k).
\end{aligned}\label{relation}
\end{equation}

By the partial subdifferential formula (see e.g.  by Proposition 2 in \cite{DCbilevel}) and \cref{g3},  there exists $\eta_{k+1} \in [0,1]$ such that
\begin{eqnarray*}
	&& \partial E_{\alpha_k}(z^{k+1},z^{k},\gamma^{k}) \\
	&=&  \partial_z E_{\alpha_k}(z^{k+1},z^{k},\gamma^{k}) \times \partial_{z_0} E_{\alpha_k}
	(z^{k+1},z^{k},\gamma^{k})\\
	&& \times \partial_\gamma E_{\alpha_k}(z^{k+1},z^{k},\gamma^{k}) \\
	&=&   \partial_z E_{\alpha_k}(z^{k+1},z^{k},\gamma^{k}) \times 
	\{ \frac{\rho}{2} (z^k- z^{k+1})\} \\
	&& \times \partial_\gamma\tilde g(z^{k+1},\gamma^{k}) \\
	&\supseteq &   \partial_z E_{\alpha_k}(z^{k+1},z^{k},\gamma^{k}) \times 
	\{ \frac{\rho}{2} (z^k- z^{k+1})\} \\
	&&  \times \{ \alpha_k \eta_{k+1}(r^{k+1} -\partial v^*(-\gamma^k)) \}.
\end{eqnarray*}

Since $z^{k+1}$ is an approximate solution to problem \cref{p-singlelevel_problem} satisfying inexact criteria \cref{p-inexact2}, there exists a vector $e_k$ such that $$e_k \in \partial {\phi}_k(z^{k+1}) + \mathcal{N}_\Sigma(z^{k+1})$$ satisfying $$\|e_k\| \le \frac{\sqrt{2}}{2} \rho\|z^k - z^{k-1}\|.$$
It follows from \cref{relation} that 
$$e_k - \frac{\rho}{2}(z^{k+1}-z^k)  \in \partial_z E_{\alpha_k}(z^{k+1},z^{k},\gamma^{k}).$$ Next,  using the fact that $-\gamma^k \in \partial v(r^{k})$ which is equivalent to $r^k\in \partial v^*(-\gamma^k)$, we have
\[ 
\begin{aligned}
	&\partial E_{\alpha_k}(z^{k+1},z^{k},\gamma^{k}) 
	\ni \begin{pmatrix}
		e_k - \frac{\rho}{2}(z^{k+1}-z^k)  \\
		\frac{\rho}{2}(z^k- z^{k+1}) \\
		\alpha_k\eta_{k+1}(r^{k+1} -  r^{k})
	\end{pmatrix}.
\end{aligned}
\]
\cref{p-re_err_eq} then  follows immediately.

\subsection{Proof of \cref{p-con_KL_alg1}}

By the assumption that the adaptive penalty sequence $\{\alpha_k\}$ is bounded, $\alpha_k = \bar{\alpha}$ for all sufficiently large $k$ and we can assume without loss of generality that 
$\alpha_k = \bar{\alpha}$ for all $k$.
Since  $L(x)$ is assumed to be bounded below,   $E_{\bar{\alpha}}(z,z_0,\gamma)$ is also bounded below. Then, according to \cref{p-suff_decreasenew}, we have  that $\lim_{k \rightarrow \infty} \|z^{k+1} - z^k\|^2 = 0$.

Let $C$ denote the set of all limit points of the  sequence $\{(z^k,z^{k-1},\gamma^{k-1})\}$. Then we have that $C$ is a closed set and $\lim_{k \rightarrow \infty} \mathrm{dist}((z^k,z^{k-1},\gamma^{k-1}), C) = 0$. And by assumption that there exists $\delta > 0$ such that $r^k_i \ge \delta$ for all $k$ and all $i=1,\dots, J$ and \cref{p-thm_limit}, any accumulation point of the sequence $\{z^k\}$ corresponds to those in the set $C$ is a KKT stationary point of problem \cref{p-singlelevel_problem}.
According to \cref{assum1}, the value function is locally Lipschitz continuous around any point in $\mathbb{R}_{> 0}^J$. Then by the assumptions that sequence $\{z^k\}$ is bounded, $r^k_i \ge \delta$ for some $\delta > 0$ and $-\gamma^{k-1} \in \partial v(z^{k-1})$, we get the boundedness of the sequence $\{\gamma^k\}$  and thus $C$ is a compact set.

Since by \cref{p-suff_decreasenew}, $E_{\bar{\alpha}}(z^k,z^{k-1},\gamma^{k-1})$ is decreasing, $E_{\bar{\alpha}}(z,z_0,\gamma)$ is bounded below and continuous on $\Sigma$ and $\lim_{k \rightarrow \infty} \|z^{k+1} - z^k\| = 0$, we have that for any subsequence $\{(z^l,z^{l-1},\gamma^{l-1})\}$ of the sequence $\{(z^k,z^{k-1},\gamma^{k-1})\}$,
\[
\bar{E} = \lim_{l \rightarrow \infty} E_{\bar{\alpha}}(z^l,z^{l-1},\gamma^{l-1}) = \lim_{k \rightarrow \infty} E_{\bar{\alpha}}(z^k,z^{k-1},\gamma^{k-1}),
\]
and thus $E_{\bar{\alpha}}$ is constant on $C$.
We can assume that $E_{\bar{\alpha}}(z^k,z^{k-1},\gamma^{k-1}) > \bar{E}$ for all $k$. Otherwise, if there exists $k > 0$ such that $E_{\bar{\alpha}}(z^k,z^{k-1},\gamma^{k-1}) = \bar{E}$, then by the above equation and \cref{p-suff_decreasenew}, we get $\|z^{k+1} - z^k\| = 0$ when $k$ is sufficiently large, which implies the convergence of sequence $\{z^k\}$ and the conclusion follows immediately. 

Since $\lim_{k \rightarrow \infty} \mathrm{dist}((z^k,z^{k-1},\gamma^{k-1}), C) = 0$ and $\lim_{k \rightarrow \infty} E_{\bar{\alpha}}(z^k,z^{k-1},\zeta^{k-1}) = \bar{E}$, for any $\epsilon, \eta > 0$, there exists $k_0$ such that $\mathrm{dist}((z^k,z^{k-1},\zeta^{k-1}), C) < \epsilon$ and $\bar{E} < E_{\bar{\alpha}}(z^k,z^{k-1},\zeta^{k-1}) < \bar{E} + \eta$ for $k \ge k_0$.
Since $E_{\bar{\alpha}}(z,z_0,\gamma)$ satisfies the Kurdyka-\L{}ojasiewicz property at each point in $C$, and $E_{\bar{\alpha}}$ is a finite constant on $C$, we can apply \cref{p-uniformKL} to obtain a continuous concave function $\varphi$ such that for any $k \ge k_0$,
\[
\begin{aligned}
	&\varphi'(E_{\bar{\alpha}}(z^k,z^{k-1},\gamma^{k-1}) - \bar{E}) \mathrm{dist} \left( 0, \partial E_{\bar{\alpha}}(z^k,z^{k-1},\gamma^{k-1}) \right) \\
	&\ge 1.
\end{aligned}
\]
Combining with \cref{p-re_err_eq} yields
\[
\begin{aligned}
	&\varphi' \left (E_{\bar{\alpha}}(z^k,z^{k-1},\gamma^{k-1}) - \bar{E}\right ) \\
	\cdot\, & \left(\frac{\sqrt{2}}{2} \rho\|z^{k-1} - z^{k-2}\| +  (\bar{\alpha}+\rho) \|z^{k} - z^{k-1}\| \right)\ge 1.
\end{aligned}
\]
The concavity of $\varphi$ and \cref{p-suff_decrease_eq2} implies 
\[
\begin{aligned}
	&\varphi'( E_{\bar{\alpha}}(z^k,z^{k-1},\gamma^{k-1}) - \bar{E}) \cdot \frac{\rho}{4} \|z^{k+1} - z^k\|^2 \\
	\le \ &\varphi'( E_{\bar{\alpha}}(z^k,z^{k-1},\gamma^{k-1}) - \bar{E})\\
	&\cdot \left( E_{\bar{\alpha}}(z^k,z^{k-1},\gamma^{k-1})  - E_{\bar{\alpha}}(z^{k+1},z^{k},\gamma^{k}) \right) \\
	\le \ &\varphi\left( E_{\bar{\alpha}}(z^k,z^{k-1},\gamma^{k-1}) - \bar{E} \right) \\
	& -  \varphi\left( E_{\bar{\alpha}}(z^{k+1},z^{k},\gamma^{k}) - \bar{E} \right).
\end{aligned}
\]
Combining the above two inequalities, we obtain 
\[
\begin{aligned}
	&\frac{4}{\rho}\left(\frac{\sqrt{2}}{2} \rho\|z^{k-1} - z^{k-2}\| +  (\bar{\alpha}+\rho) \|z^{k} - z^{k-1}\| \right) \\
	&\cdot \Bigg[\varphi\left( E_{\bar{\alpha}}(z^k,z^{k-1},\gamma^{k-1}) - \bar{E} \right) \\
	&\hspace{70pt}-  \varphi\left( E_{\bar{\alpha}}(z^{k+1},z^{k},\gamma^{k}) - \bar{E} \right) \Bigg] \\
	\ge\, &\|z^{k+1} - z^k\|^2.
\end{aligned}
\]

Multiplying both sides of this inequality by 4 and taking the square root, and by $2ab \le a^2 + b^2$, we have
\[
\begin{aligned}
	4\|z^{k+1} - z^k\| \le \,&\|z^{k} - z^{k-1}\| + \|z^{k-1} - z^{k-2}\|\\  &\hspace{-20pt} + \frac{16(\bar{\alpha}+\rho)}{\rho}\Big[\varphi\left( E_{\bar{\alpha}}(z^k,z^{k-1},\gamma^{k-1}) - \bar{E} \right) \\
	&  -  \varphi\left( E_{\bar{\alpha}}(z^{k+1},z^{k},\gamma^{k}) - \bar{E} \right) \Big].
\end{aligned}
\]
Summing up the above inequality for $i = k_0 + 1, \ldots, k$, we have
\[
\begin{aligned}
	\sum_{i = k_0}^k 4\|z^{i+1} - z^i\| \le  &	\sum_{i = k_0}^k \left( \|z^{i} - z^{i-1}\| + \|z^{i-1} - z^{i-2}\| \right)\\
	&\hspace{-45pt}+ \frac{16(\bar{\alpha}+\rho)}{\rho} \Big[\varphi\left( E_{\bar{\alpha}}(z^{k_0},z^{k_0-1},\gamma^{k_0-1}) - \bar{E} \right) \\
	& \hspace{20pt}- \varphi\left( E_{\bar{\alpha}}(z^{k+1},z^{k},\gamma^{k}) - \bar{E} \right) \Big],
\end{aligned}
\]
and since $\varphi \ge 0$, we get
\[
\begin{aligned}
	\sum_{i = k_0}^k 2\|z^{i+1} - z^i\| \le  &2\|z^{k_0} - z^{k_0-1}\| + \|z^{k_0-1} - z^{k_0-2}\| \\ &\hspace{-45pt}+ \frac{16(\bar{\alpha}+\rho)}{\rho}\varphi\left( E_{\bar{\alpha}}(z^{k_0},z^{k_0-1},\gamma^{k_0-1}) - \bar{E} \right).
\end{aligned}
\]
Taking $k \rightarrow \infty$ in the above inequality shows that
\[
\sum_{k = 1}^\infty \|z^{k+1} - z^k\| < \infty,
\]
and thus the sequence $\{z^k\}$ is a Cauchy sequence. Hence the sequence $\{z^k\}$ is convergent and we get the conclusion.

\section{Detailed description of  VF-iDCA in the general setting}\label{sectionB}

Our value function based DC algorithm as well as its convergence analysis can be straightforwardly extended to BLP in a more general setting with LL constraints.
\begin{equation}\label{general_problem}
	\begin{aligned}
		\min_{x \in X, u \in U, \lambda \in \mathbb{R}_+^J} ~~& L(x,u) \\
		\text{s.t.} ~~& x \in \underset{x' \in X}{\mathrm{argmin}} \bigg\{ l(x',u) + \sum_{i = 1}^{J}\lambda_i P_i(x',u) \\
		& \hspace{35pt} \text{s.t.} \quad  g(x', u) \le 0 \bigg\},
	\end{aligned}
\end{equation}
where $X \subseteq \mathbb{R}^n$, $U \subseteq \mathbb{R}^d$ are closed convex sets, $L, l : \mathbb{R}^n \times \mathbb{R}^d \rightarrow \mathbb{R}$, 
$P_i : \mathbb{R}^n \times \mathbb{R}^d \rightarrow \mathbb{R}_+,~ i = 1,\ldots,J$, $g_i : \mathbb{R}^n \times \mathbb{R}^d \rightarrow \mathbb{R},~ i = 1,\ldots,m$ are  convex functions defined on an open convex set containing  $X\times U$ and $g: \mathbb{R}^n \times \mathbb{R}^d \rightarrow \mathbb{R}^m $ is defined as $g(x,u) = (g_1(x,u),\ldots,g_m(x,u))$.  To ensure the convergence, we assume that 
function $L(x,u)$ is bounded below on an open convex set containing  $X\times U$.

Decoupling the hyperparameter variable $\lambda$ from the regularization term, and introducing a new variable $r$ result in the following BLP:
\begin{equation}\label{general_problem_c}
	\begin{aligned}
		\min_{x \in X, u \in U, \lambda \in \mathbb{R}_+^J} ~~& L(x,u) \\
		\text{s.t.} ~~& x \in \underset{x' \in X}{\mathrm{argmin}} \bigg\{ l(x',u)  \\
		& \hspace{35pt} \text{s.t.} \quad  g(x', u) \le 0, \\
		& P_i(x',u) \le r_i,~ i=1,\ldots,J\bigg \}.
	\end{aligned}
\end{equation}
Denote the solution sets of the lower level programs in problem \cref{general_problem} and \cref{general_problem_c} by $S_p(\lambda, u) $ and $S_c(r, u) $, respectively.

BLP \cref{general_problem_c} is then  reformulated as  the following single-level DC program:
\begin{equation}\label{singlelevel_problem}
	\begin{aligned}
		\min_{x \in X, u \in U, \lambda \in \mathbb{R}_+^J} ~~& L(x,u) \\
		\text{s.t.} ~~~& l(x,u) - v(r,u) \le 0, g(x, u) \le 0 ,\\
		&P_i(x,u) \le r_i,~ i=1,\ldots,J,\\
	\end{aligned}
\end{equation}
where $v(r,u)$ is the value function of the LL problem governed by $(r,u)$ $$
\begin{aligned}
	&v(r,u) := \min_{x' \in X} ~\{l(x',u)~\text{s.t.}~ g(x', u) \le 0 , \\
	& \hspace{100pt} P_i(x',u) \le r_i,~ i=1,\ldots,J \}.
\end{aligned}
$$ 

To ensure the convergence,  assume that for each $r$ belonging to any open subset of $ \mathbb{R}_+^J$,   the feasible region ${\cal F}(r,u):=\{ x \in X  \text{s.t.} g(x,u)\leq 0, P_i(x,u)\leq r_i, i=1,\dots, J\} \not =\emptyset$ and $l(x,u)$ is bounded on ${\cal F}(r,u)$. 
Thanks to the full convexity after hyperparameter variable decoupling, $v(r,u)$ is indeed convex and locally Lipschitz continuous around every point in set $\mathbb{R}^J_{>0} \times U$ (see, e.g., Lemma 3 in \cite{DCbilevel}). 

Given a current iteration $({x}^k,u^k,r^k)$ for each $k$, solving the LL problem parameterized by $u^k$ and $r^k$
\begin{equation}\label{LL-iter}
	\begin{aligned}
	\min_{x \in X} ~~ l(x,u^k)~\text{s.t.}~ &g(x, u^k) \le 0 , \\
	&  P_i(x,u^k) \le r_i^k,~ i=1,\ldots,J,
	\end{aligned}
\end{equation}
leads to a solution $\tilde{x}^k \in S_c(r^k, u^k)$ and a corresponding KKT multiplier $(\gamma^k,\zeta^k) \in \mathcal{M}(\tilde{x}^k,u^k,r^k)$, where $\mathcal{M}(x,u,r)$ denotes the set of multipliers of the lower level problem, 
\[
\begin{aligned}
	&\mathcal{M}(x,u,r):=\Big \{ \lambda\in \mathbb{R}^J_+, \mu \in \mathbb{R}^m_+ \mid 0\in \partial_x l(x,u) + \mathcal{N}_X(x) \\
	&+\sum_{i=1}^J \lambda_i \partial_x P_i(x,u) + \sum_{i=1}^m\mu_i\partial_x g_i(x,u), ~ \langle \mu, g(x,u) \rangle  = 0,\\
	&\hspace{85pt} ~ \lambda_i (P_i(x,u)-r_i)=0, i=1,\dots, J \Big \}.
\end{aligned}
\]
Select 
\begin{equation}
	\begin{aligned}
		\xi^k \in  \partial_u l(\tilde{x}^k,u^k) &+ \sum_{i = 1}^{J}\gamma_i ^k\partial_u P_i(\tilde{x}^k,u^k) \\ &\hspace{50pt} + \sum_{i=1}^m\zeta^k_i\partial_u g_i(\tilde{x}^k,u^k).
		 \label{subgvaluefunction}
	\end{aligned}
\end{equation}
Then by sensitivity analysis (see, e.g., Theorem 3 in \cite{DCbilevel}), it can be easily checked that $(- \gamma^k, \xi^k) \in \partial v(r^k,u^k)$.
Compute $z^{k+1} :=(x^{k+1},u^{k+1},r^{k+1})$ as an approximate minimizer of the strongly convex subproblem			
\begin{equation} \label{DCA2_subproblem}
	\begin{aligned}
		\min_{ \substack{x \in X, u \in U,r \in \mathbb{R}_+^J}} ~~ {\phi}_k(x,u,r) := &L(x,u) +  \frac{\rho}{2} \| z - z^k \|^2 \\
		& \hspace{-100pt} +\alpha_k \max_{ \substack{j=1,\ldots,J\\ i =1,\ldots,m}}\{0, {V}_k(x,u,r), P_j(x,u) - r_j, g_i(x,u)\},
	\end{aligned}
\end{equation}
where $z:= (x,u,r)$, $z^k := (x^k,u^k,r^k)$ and $${V}_k(x,u,r):= l(x,u) - l(\tilde{x}^k,u^k) - \langle \xi^{k}, u - u^k\rangle +\langle \gamma^k, r-r^k\rangle. $$
Denoting $\Sigma:=  X\times  U\times \mathbb{R}_+^J $, we may introduce an inexact condition for choosing $z^{k+1}$
\begin{equation}\label{inexact2}
	\mathrm{dist}(0, \partial \phi_k(z^{k+1}) + \mathcal{N}_\Sigma(z^{k+1})) \le \frac{\sqrt{2}}{2} \rho\|z^k - z^{k-1}\|.
\end{equation}

Using above constructions, 
letting $$
\begin{aligned}
	t^{k+1}
	:= \max_{ \substack{j=1,\ldots,J\\ i =1,\ldots,m}} \Big\{0, &{V}_k(z^{k+1},u^{k+1},r^{k+1}),  \\ &\hspace{-25pt}P_j(x^{k+1},u^{k+1}) - r^{k+1}_j,  g_i(x^{k+1},u^{k+1}) \Big\},
\end{aligned}$$
we are ready to present  VF-iDCA in \cref{ipDCA}.
\begin{algorithm}[h]
	\caption{VF-iDCA}\label{ipDCA}
	\begin{algorithmic}[1]
		\STATE Take an initial point $(x^0,u^0,r^0) \in X \times U \times \mathbb{R}_+^J$; $c_\alpha, \delta_\alpha > 0$; an initial penalty parameter $\alpha_0 > 0$; tolerance $tol > 0$.
		\FOR{$k=0,1,\ldots$}
		\STATE Solve LL problem \cref{LL-iter}. Find
		$\tilde{x}^k \in S_c(r^k,u^k)$ and KKT multiplier $\gamma^k,\zeta^k$.		
		\STATE Solve problem \cref{DCA2_subproblem} up to tolerance in \cref{inexact2}. Find an approximate solution $(x^{k+1}, u^{k+1}, r^{k+1})$.
		\STATE Stopping test. 
		
		Stop if $\max\{ \|z^{k+1} - z^k\|,t^{k+1} \} < tol$.
		\STATE Adaptive penalty parameter update.
		Set
		\begin{equation*}
			\alpha_{k+1} = \left\{
			\begin{aligned}
				&\alpha_k + \delta_\alpha, &&\text{if}~\max\{\alpha_k, 1/t^{k+1}\} < \frac{c_\alpha}{\Delta^{k+1}}, \\
				&\alpha_k, &&\text{otherwise},
			\end{aligned}\right.
		\end{equation*}
		where $\Delta^{k+1}:=\|z^{k+1} - z^k\|$.
		\ENDFOR
	\end{algorithmic}
\end{algorithm}

Obviously  similar convergence results for the general problem as in  Theorems \ref{p-thm_limit} and \ref{p-con_KL_alg1} can  be obtained.
\section{Hyperparameter Decoupling Illustration}\label{sectionC}

We illustrate  how to adopt the hyperparameter decoupling technique and cast the BLPs listed in Table \ref{table:model} in the formulation of  \cref{reform_problem} with fully convex LL problem.

\paragraph{Elastic net \cite{zou2003regression}} is a regularized linear regression method that combines the ridge and lasso penalties. The hyperparameter optimization problem for elastic net can be formulated as the following BLP:
\begin{align}
\min_{\bm{\beta} \in \mathbb{R}^p, \lambda\in \mathbb{R}_+^2}\ & \frac12 \sum_{i \in I_{\text{val}}} | b_i - \bm{\beta}^\top \mathbf{a}_i |^2 \\
\text{s.t. } \bm{\beta} \in \mathop{\arg\min}_{\hat{\bm{\beta}} \in \mathbb{R}^p} & \bigg\{ 
\frac12 \sum_{i\in I_{\text{tr}}} | b_i - \hat{\bm{\beta}}^\top \mathbf{a}_i |^2 \nonumber \\
& \qquad + \lambda_1\|\hat{\bm{\beta}}\|_1 + \frac{\lambda_2}2 \|\hat{\bm{\beta}}\|_2^2 \bigg\}. \nonumber
\end{align}

For the implementation of VF-iDCA, we adopt the hyperparameter decoupling and thus work on the following BLP:
\begin{align}
    \min_{\bm{\beta} \in \mathbb{R}^p, \mathbf{r} \in \mathbb{R}_+^2}\ & \frac12 \sum_{i \in I_{\text{val}}} | b_i - \bm{\beta}^\top \mathbf{a}_i |^2 \\
    \text{s.t. } \bm{\beta} \in & \mathop{\arg\min}_{\hat{\bm{\beta}} \in \mathbb{R}^p} \bigg\{ 
    \frac12 \sum_{i\in I_{\text{tr}}} | b_i - \hat{\bm{\beta}}^\top \mathbf{a}_i |^2 \nonumber \\ 
    & \qquad \quad \text{ s.t. }
    \|\hat{\bm{\beta}}\|_1 \le r_1,\ \frac12 \|\hat{\bm{\beta}}\|_2^2 \le r_2 \bigg\}. \nonumber
\end{align}

\paragraph{Sparse group lasso \cite{simon2013sparse}} is a regularized linear regression method well-suited for cases where the features have a natural grouping. The sparse group lasso penalty involves both $\|\cdot\|_2$ and $\|\cdot\|_1$ norms, and encourages sparsity both on the group level and on the within group level.

Suppose the $p$ features are divided into $M$ groups. Denoting the corresponding coefficients of $m$-th group as $\bm{\beta}^{(m)}$ and the corresponding regularization parameter as $\lambda_m$ ($r_m$ for the VF-iDCA), we can write the hyperparameter optimization problem  in the following way:

\begin{align}
    \min_{\bm{\beta} \in \mathbb{R}^p, \bm{\lambda}\in \mathbb{R}_+^{M+1}}\ & \frac12 \sum_{i \in I_{\text{val}}} | b_i - \bm{\beta}^\top \mathbf{a}_i |^2 \\
    \text{s.t. }\ \bm{\beta} \in & \mathop{\arg\min}_{\hat{\bm{\beta}} \in \mathbb{R}^p} \bigg\{ \frac12 \sum_{i\in I_{\text{tr}}} | b_i - \hat{\bm{\beta}}^\top \mathbf{a}_i |^2 \nonumber \\ 
    & \qquad + \sum_{m=1}^M \lambda_m \|\hat{\bm{\beta}}^{(m)}\|_2 + \lambda_{M+1}\|\hat{\bm{\beta}}\|_1 \bigg\}. \nonumber
\end{align}
Note that in the above, we consider a variation on \cite{simon2013sparse} where we use a separate regularization parameter for each of the groups, rather than a single regularization parameter for the sum of the penalties across all groups, along the lines of \cite{feng2018gradient}. 

For the implementation of VF-iDCA, we adopt the hyperparameter decoupling and thus work on the following BLP:
\begin{align}
    \min_{\bm{\beta} \in \mathbb{R}^p, \mathbf{r} \in \mathbb{R}_+^{M+1}} & \frac12 \sum_{i \in I_{\text{val}}} \| b_i - \bm{\beta}^\top \mathbf{a}_i \|^2 \\
    \text{s.t. }\ \bm{\beta} \in \mathop{\arg\min}_{\hat{\bm{\beta}} \in \mathbb{R}^p} & \bigg\{ \frac12 \sum_{i\in I_{\text{tr}}} \| b_i - \hat{\bm{\beta}}^\top \mathbf{a}_i \|^2 \nonumber \\ 
    \text{s.t. } \|\hat{\bm{\beta}}^{(m)}\|_2 & \le r_m,\ m = 1,\dots,M,\ \|\hat{\bm{\beta}}\|_1 \le r_{M+1} \bigg\}. \nonumber
\end{align}

\paragraph{Low-rank matrix completion} is a flexible framework for reduced-rank modelling of matrix-valued data. We consider a variation explored in \citet{feng2018gradient}, where we have access to additional information corresponding to each row and column. Specifically, suppose that for matrix $M\in \mathbb{R}^{n\times n}$, we observe some entries $M_{ij}$ where $(i, j) \in \Omega$ and do not have access to the rest. We denote the row features $X\in \mathbb{R}^{n\times p}$ and column features as $Z\in \mathbb{R}^{n\times p}$. We model the matrix as the sum of a low rank effect $\Gamma$ and a linear combination of the row features and the column features. Denote the coefficients of row features and column features by $\bm{\theta}$ and $\bm{\beta}$, respectively. Furthermore, suppose that there is a natural grouping of the row and column features, which partitions the corresponding coefficients as $\{\bm{\theta}^{(g)}\}_{g=1}^{G}$ and $\{\bm{\beta}^{(g)}\}_{g=1}^{G}$. So with the nuclear norm penalty on $\Gamma$ and the group lasso penalty on the linear part, we have the following model:
\begin{align*}
	&\min_{\substack{\bm{\theta}\in \mathbb{R}^p,\ \bm{\beta}\in \mathbb{R}^p,\\ \Gamma\in \mathbb{R}^{n\times n},\ \bm{\lambda}\in \mathbb{R}^{2G+1}_+}}
	\sum_{(i, j) \in \Omega_{\text{val}}} \frac12|M_{ij} - \mathbf{x}_i \bm{\theta} - \mathbf{z}_j \bm{\beta} - \Gamma_{ij} |^2\\
	&\text{s.t. } \bm{\theta},\ \bm{\beta},\ \Gamma \in \\
	& \mathop{\arg\min}_{\substack{\hat{\bm{\theta}}\in \mathbb{R}^p,\ \hat{\bm{\beta}}\in \mathbb{R}^p,\\ \hat{\Gamma}\in \mathbb{R}^{n\times n}}} 
	\bigg\{ \sum_{(i, j) \in \Omega_{\text{tr}}} \frac12|M_{ij} - \mathbf{x}_i \hat{\bm{\theta}} - \mathbf{z}_j \hat{\bm{\beta}} - \hat{\Gamma}_{ij} |^2 \\
	& \hspace{20pt} + \lambda_0 \|\hat{\Gamma}\|_* + \sum_{g=1}^G \lambda_g \|\hat{\bm{\theta}}^{(g)}\|_2 + \sum_{g=1}^G \lambda_{g+G} \|\hat{\bm{\beta}}^{(g)}\|_2 \bigg\},
\end{align*}
where $\Omega = \Omega_{\text{val}} \cup \Omega_{\text{tr}}$, $\Omega_{\text{val}} \cap \Omega_{\text{tr}} = \emptyset$.

For the implementation of VF-iDCA, we adopt the hyperparameter decoupling first and obtain the following BLP:
\begin{align*}
	&\min_{\substack{\bm{\theta}\in \mathbb{R}^p,\ \bm{\beta}\in \mathbb{R}^p,\\ \Gamma\in \mathbb{R}^{n\times n},\ \mathrm{r}\in \mathbb{R}^{2G+1}_+}}
	 \sum_{(i, j) \in \Omega_{\text{val}}} \frac12|M_{ij} - \mathbf{x}_i \bm{\theta} - \mathbf{z}_j \bm{\beta} - \Gamma_{ij} |^2\\
	&\text{s.t. } \bm{\theta}, \bm{\beta}, \Gamma \in \\ &\mathop{\arg\min}_{\substack{\hat{\bm{\theta}}\in \mathbb{R}^p, \hat{\bm{\beta}}\in \mathbb{R}^p,\\ \hat{\Gamma}\in \mathbb{R}^{n\times n}}} 
	\bigg\{ \sum_{(i, j) \in \Omega_{\text{tr}}} \frac12|M_{ij} - \mathbf{x}_i \hat{\bm{\theta}} - \mathbf{z}_j \hat{\bm{\beta}} - \hat{\Gamma}_{ij} |^2 \\
	& \hspace{40pt}\text{s.t. } \|\hat{\Gamma}\|_* \le r_0,\ \|\hat{\bm{\theta}}^{(g)}\|_2 \le r_g,\ g = 1,\dots, G\\ 
	& \hspace{80pt} \|\hat{\bm{\beta}}^{(g)}\|_2 \le r_{g+G},\ g = 1,\dots, G \bigg\}
\end{align*}

\paragraph{Support vector machines \cite{cortes1995support}} are supervised learning methods for binary classification. Support vector machines solve
\[
\min_{\substack{-\bar{\mathbf{w}} \le \mathbf{w} \le \bar{\mathbf{w}} \\ c \in \mathbb{R}}}
\bigg\{ \frac{\lambda}{2}\|\mathbf{w}\|^2 + \sum_{j\in \Omega_{\text{tr}}}\max( 1 - b_j(\mathbf{w}^\top \mathbf{a}_j - c),0) \bigg\}.    
\]
We note that this is a variant of support vector machines proposed by   \citealt{kunapuli2008bilevel} which incorporates feature selection through the box constraints on $\mathbf{w}$. 

Recent work \cite{kunapuli2008classification} considered the $T$-fold cross validation method for selecting the hyperparameters $\lambda$ and $\bar{\mathbf{w}}$, and treated it as a bilevel optimization problem. 
A given data set $\Omega$ is randomly partitioned into $T$ pairwise disjoint subset called the validation sets $\{\Omega_{\text{val}}^t\}_{t=1}^T$. For each validation set $\Omega_{\text{val}}^t$, the corresponding training set is $\Omega_{\text{tr}}^t := \Omega \backslash \Omega_{\text{val}}^t$. We select the best hyperparameters via minimizing the hinge loss based on the validation set and its partition by following BLP:
\begin{equation*}
	\begin{aligned}
		&\min
		_{\lambda,\bar{\mathbf{w}},\mathbf{w}^1,\dots, \mathbf{w}^T, \mathbf{c}} 
		~~ 
		\Theta( \mathbf{w}^1, \dots,\mathbf{w}^T ,\mathbf{c})\\
		&~~ \qquad \begin{aligned}
			\text{s.t.} ~~ & \lambda \ge 0, \quad \bar{\mathbf{w}}_{lb}\le \bar{\mathbf{w}} \le \bar{\mathbf{w}}_{ub},   \\
			& \mathrm{and~for}~t = 1,\ldots,T:\\
			& (\mathbf{w}^t, c^t) \in \underset{\tiny \begin{matrix}
					-\bar{\mathbf{w}} \le \mathbf{w} \le \bar{\mathbf{w}} \\ c \in \mathbb{R}	\end{matrix} }{\mathrm{argmin}} \Bigg\{ \frac{\lambda}{2}\|\mathbf{w}\|^2 \\
			&\hspace{30pt}+ \sum_{j\in \Omega_{trn}^t}\max( 1 -b_j( \mathbf{a}_j^T\mathbf{w}-c),0) \Bigg\},
		\end{aligned}
	\end{aligned}
\end{equation*}
with 
\[
\begin{aligned}
	&\Theta( \mathbf{w}^1, \dots,\mathbf{w}^T ,\mathbf{c}) \\
	:= \,&\frac{1}{T} \sum_{t=1}^{T} \frac{1}{|\Omega_{val}^t|} \sum_{j\in \Omega_{val}^t}\max(1-b_j(\mathbf{a}_j^T\mathbf{w}^t
	-c^t
	),0), 
\end{aligned} 
\]
where $|\Omega|$ denotes the number of elements in set $\Omega$. 

For the implementation of VF-iDCA, we adopt the hyperparameter decoupling first and obtain the following BLP:
\begin{align*}
	& \min
	_{r,\bar{\mathbf{w}},\mathbf{w}^1,\dots, \mathbf{w}^T, \mathbf{c}} 
	~~ 
	\Theta( \mathbf{w}^1, \dots,\mathbf{w}^T ,\mathbf{c})\\
	& \begin{aligned}
		\text{s.t.} ~~ & r \ge 0, \quad \bar{\mathbf{w}}_{lb}\le \bar{\mathbf{w}} \le \bar{\mathbf{w}}_{ub},   \\
		& (\mathbf{w}^t, c^t) \in \\
		& \mathop{\arg\min}_{\mathbf{w} \in \mathbb{R}^p, c \in \mathbb{R}}
		\Bigg\{ \sum_{j\in \Omega_{trn}^t}\max( 1 -b_j( \mathbf{a}_j^T\mathbf{w}-c),0) \\
		& \hspace{7em} \text{s.t. } \frac{1}{2}\|\mathbf{w}\|^2\le r, 
		-\bar{\mathbf{w}} \le \mathbf{w} \le \bar{\mathbf{w}} \Bigg\}.
	\end{aligned}
\end{align*}

\section{Numerical experiments}\label{sectionD}

\subsection{Numerical experiments on synthetic data}

\subsubsection{Elastic Net}

We simulate data in a similar manner as \citet{feng2018gradient} as follows.
We draw $\mathbf{a}_i\in \mathbb{R}^{p}$ from a $N(\bm{0}, \bm{I})$ distribution with $\mathrm{cor}(a_{ij}, a_{ik}) = 0.5^{|j - k|}$. 
We draw the response $\mathbf{b}$ according to $b_i = \bm{\beta}^\top \mathbf{a}_i + \sigma \epsilon_i$, where $\bm{\beta}$ is randomly generated such that $\beta_i$ is either 0 or 1 and $\sum_{i = 1}^{p}\beta_i = 15$; $\bm{\epsilon}$ is sampled from the standard Gaussian distribution, and $\sigma$ was chosen such that the signal-to-noise ratio $\text{SNR} \overset{\Delta}{=} \|A\bm{\beta}\|/\|\mathbf{b} - A\bm{\beta}\|$ was 2.

For grid search, random search and TPE, we defined $u_1 = \log_{10}(\lambda_1)$ and $u_2 = \log_{10}(\lambda_2)$, and searched over $u_1, u_2 \in [-5, -2]$, as in \cite{feng2018gradient}. Grid search was performed on a $10\times 10$ uniformly-spaced grid. Random search was performed with 100 uniform random samples. The space used in TPE for both $u_1$ and $u_2$ was a uniform distribution on $[-5, 2]$. For IGJO, the initial guesses for $\lambda_1$ and $\lambda_2$ were $0.01$ and $0.01$, respectively. For IFDM, the inital guesses for $\lambda_1$ and $\lambda_2$ were were $0.01\lambda_{\text{max}}$ as in \cite{bertrand2021implicit}, where $\lambda_{\text{max}} = 1/|I_{\text{tr}}| \max( A_{\text{tr}}^\top \mathbf{b}_{\text{tr}} )$. The maximum number of iterations of IFDM was set to be 50. For VF-iDCA, the initial guesses for $r_1$ and $r_2$ were $10$ and $5$, respectively. VF-iDCA was stopped when 
\begin{equation}
	\label{stopping_criterion}
	\max\left\{
	\frac{\lVert \mathbf{z}^{k+1} - \mathbf{z}^k\rVert}{\sqrt{1+\|\mathbf{z}^k\|^2}},
	t^{k+1}
	\right\} < tol,
\end{equation}
and we set $tol = 0.1 $.

\subsubsection{Sparse Group Lasso}

The data generation method was mainly referred to \citet{feng2018gradient} as follows. Each dataset contains 100 training data, 100 validation data and 100 test data. Each $\mathbf{a}_i\in \mathbb{R}^{p}$ were sampled from the standard normal distribution. 
The response $\bm{b}$ was generated by $b_i = \bm{\beta}^\top \mathbf{a}_i + \sigma \epsilon_i$, where $\bm{\beta} = \big[\bm{\beta^{(1)}}, \bm{\beta^{(2)}},\bm{\beta^{(3)}}\big]$, $\bm{\beta}^{(i)} = (1,2,3,4,5,0,\dots,0)$, for $i = 1,2,3$. $\bm{\epsilon}$ are generated from the standard normal distribution, and $\sigma$ was chosen such that the SNR is 2.

For grid search, we did the search on two hyperparameters $\mu_1$, $\mu_2$ such that $\lambda_m = 10^{\mu_1}$, for $m = 1,2,\dots,M$ and $\lambda_{M+1} = 10^{\mu_2}$, we use a 10$\times$10 uniform-spaced grid on $[-3, 1]\times [-3,1]$. For random search and TPE method, we search over the transformed variables $u_m$, where $u_m = \log_{10}(\lambda_m)$, for $m = 1,2,\dots, M+1$, and the space of $u_m$ is defined as a uniform distribution on $[-3, 1]$. For implicit differentiation method, we use the initial guess $[0.01, 0.01,\dots, 0.01]$ for $\bm{\lambda}$. For VF-iDCA, we use the initial guess $[10,10,\dots,10]$ for $\mathbf{r}$.
The stopping criterion used for VF-iDCA here was \cref{stopping_criterion} with $tol = 0.05$.
In all experiments, the features are grouped by order, i.e., the first $p/M$ features are grouped as the first group, the next $p/M$ features are grouped as the second group, etc.

\subsubsection{Low-rank matrix completion}
The data generation method was also referred to \citet{feng2018gradient} as follows. We took two entries per row and column as the training set $\Omega_{\text{tr}}$ and one entry per row and column as the validation set $\Omega_{\text{val}}$. The rest entries were collected as test set $\Omega_{\text{test}}$. The row features are grouped into 12 groups of 3 covariates each, and same for the columns features, i.e., $p = 36,\ G = 12$.

The true coefficients chosen as $\bm{\alpha}^{(g)} = g \mathbf{1}_3$ for $g = 1,\dots,4$ and $\bm{\beta}^{(g)} = g \mathbf{1}_3$ for $g = 1,\ 2$. The rest of the coefficients were zero. We generated rank-one effect matrices $\Gamma = \mathbf{u} \mathbf{v}^\top$, where $\mathbf{u}$ and $\mathbf{v}$ were sampled from a standard normal distribution. The row features and column features $X$ and $Z$ were sampled from a standard normal distribution, and scaled such that the Frobenius norm of $X \bm{\alpha} \mathbf{1}^\top + (Z \bm{\beta} \mathbf{1}^\top)^\top$ was the same as $\Gamma$. The matrix data were generated by $M_{ij} = \mathbf{x}_i \bm{\alpha} + \mathbf{z}_j \bm{\beta} + \Gamma_{ij} + \sigma \epsilon_{ij}$, where $\epsilon_{ij}$ were generated from the standard normal distribution, and $\sigma$ was chosen such that the SNR is 2.

For grid search, we did the search on two hyperparameters $\mu_1$, $\mu_2$ such that $\lambda_0 = 10^{\mu_1}$ and $\lambda_g = 10^{\mu_2}$, for $g = 1,\dots, 2G$, we use a $10\times 10$ uniform-spaced grid on $[-3.5, -1] \times [-3.5, -1]$ as \citet{feng2018gradient} did. For random search and TPE method, we search over the transformed variables $u_g$, where $u_g = \log_{10}(\lambda_m)$, for $m = 0,1,2,\dots, 2G$, and the space of $u_g$ is defined as a uniform distribution on $[-3.5, -1]$. For IGJO, we use the initial guess $[0.005, 0.005,\dots, 0.005]$ for $\bm{\lambda}$. For VF-iDCA, we use the initial guess $[1, 0.1, 0.1,\dots,0.1]$ for $\mathbf{r}$. The stopping criterion used for VF-iDCA here was \cref{stopping_criterion} with $tol = 0.05$. In all experiments, the features are grouped by order, i.e., the first $3$ features are grouped as the first group, the next $3$ features are grouped as the second group, etc.

\subsection{Application to real data}

\subsubsection{Elastic Net}

The numerical settings for these method are kept the same as in \cref{numerical_EN}.
For grid search, random search and TPE method, we did the search on the logarithmic hyperparameters $u_1 = \log_{10}(\lambda_1)$ and $u_2 = \log_{10}(\lambda_2)$, and the range of both $u_1$ and $u_2$ is $[-5,2]$ as in \citet{feng2018gradient}. Grid search was performed on a $10\times 10$ uniform-spaced grid. Random search was performed with 100 uniformly random sampling. The space used in TPE method for both $u_1$ and $u_2$ is deifned as a uniform distribution on $[-5, 2]$. For IGJO, the initial guesses for $\lambda_1$ and $\lambda_2$ were $0.01$ and $0.01$, respectively. For IFDM, the inital guesses for $\lambda_1$ and $\lambda_2$ were $0.01\lambda_{\text{max}}$ as in \citet{bertrand2021implicit}, where $\lambda_{\text{max}} = 1/|I_{\text{tr}}| \max( A_{\text{tr}}^\top \mathbf{b}_{\text{tr}} )$. The maximum number of iterations of IFDM was set to be 10. For VF-iDCA, the initial guesses for $r_1$ and $r_2$ are 10, 5. VF-iDCA was stopped when  \cref{stopping_criterion} was satisfied with $tol = 0.1$.

\subsubsection{Support Vector Machine}

For grid search and random search, we did the search over two hyperparameters $\mu_1$ and $\mu_2$, where $\lambda = 10^{\mu_1}$, $\bar{\mathbf{w}} = (10^{\mu_2}, \dots, 10^{\mu_2})^\top$, the range of $\mu_1$ is -4 to 4, and the range of $\mu_2$ is -6 to 2. For TPE method, we search the $\log_{10}(\lambda)$ in $[-4, 4]$, and $\log_{10}(\bar{w}_i)$ in $[-6, 2]$. However, since TPE is slow when the dimension is high, we set the maximum number of iteration to be 10. And we also tested TPE method on the simplified model, i.e., the same setting as the search methods, and denote such a method as TPE2. For the simplified version, we set the maximum number of iterations to be 100. For VF-iDCA, the initial guess of $r$ is 10 and the initial guess of $\bar{\mathbf{w}}$ is $(10^{-6}, \cdots, 10^{-6})^\top$. 
We implement VF-iDCA with two different stopping criteria, i.e., \cref{stopping_criterion} with $tol = 0.01$ and $tol = 0.1$. We denote the one with $tol=0.01$ by VF-iDCA, and the one with $tol=0.1$ by VF-iDCA-t.
We set $\bar{\mathbf{w}}_{lb} = 10^{-6}$ and $\bar{\mathbf{w}}_{ub} = 10$ in the model.


\end{document}